\documentclass[a4paper,10pt]{article}

\usepackage{euscript}
\usepackage{amsfonts}
\usepackage{amsmath,amscd, amsthm, amssymb}
\usepackage[english]{babel}
\usepackage[latin1]{inputenc}
\usepackage{graphicx}
\usepackage[arrow, matrix, curve]{xy}
\usepackage{dsfont}\let\mathbb\mathds
\usepackage{rotating}
\usepackage{url}

\newtheorem{thm}{Theorem}
\newtheorem{prop}[thm]{Proposition}
\newtheorem{cor}[thm]{Corollary}
\newtheorem{lem}[thm]{Lemma}

\newtheorem{Def}[thm]{Definition}

\def\RR{{\mathbb{R}}}

\def\cartesien{%
    \ar@{-}[]+R+<6pt,-1pt>;[]+RD+<6pt,-6pt>%
    \ar@{-}[]+D+<1pt,-6pt>;[]+RD+<6pt,-6pt>%
  }
\newcommand{\hook}[1][r]
   {\ar@{}[#1] |*[o][F]{\hbox{%
         \vrule width 1.5mm height 0pt depth 0pt%
         \vrule width 0pt height .75mm depth .75mm%
         }}
     \ar@{^{(}->}[#1]}
\newcommand{\hookl}[1][r]
   {\ar@{}[#1] |*[o][F]{\hbox{%
         \vrule width 1.5mm height 0pt depth 0pt%
         \vrule width 0pt height .75mm depth .75mm%
         }}
     \ar@{_{(}->}[#1]}     
\newcommand{\demde}[1]{\begin{proof} de #1}
\newcommand{\dem}{\begin{proof}}
\newcommand{\cqfd}{\end{proof}}
\let\cat\mathfrak 
  \newcommand{\UN}[4][r]{%
    \ar@/^1pc/[#1]^{#2}_*=<0.3pt>{}="HAUT"
    \ar@/_1pc/[#1]_{#3}^*=<0.3pt>{}="BAS"
    \ar @{=>} "HAUT";"BAS" ^{#4}
  }
\newcommand{\DEUX}[6][r]{
    \ar@/^2pc/[#1]^{#2}_*=<0.3pt>{}="HAUT"
    \ar@{}    [#1]     ^*=<0.3pt>{}="MILIEUHAUT"
                       _*=<0.3pt>{}="MILIEUBAS"
    \ar[#1]_(0.3){#3}                  
    \ar@/_2pc/[#1]_{#4}^*=<0.3pt>{}="BAS"
    \ar @{=>} "HAUT";"MILIEUHAUT" ^{#5}
    \ar @{=>} "MILIEUBAS";"BAS" ^{#6}
  }   
 \newcommand{\eq}[1][r]
   {\ar@<-3pt>@{-}[#1]
    \ar@<-1pt>@{}[#1]|<{}="gauche"
    \ar@<+0pt>@{}[#1]|-{}="milieu"
    \ar@<+1pt>@{}[#1]|>{}="droite"
    \ar@/^2pt/@{-}"gauche";"milieu"
    \ar@/_2pt/@{-}"milieu";"droite"}
 \newcommand{\incl}[1][r]
  {\ar@<-0.2pc>@{^(-}[#1] \ar@<+0.2pc>@{-}[#1]}

\newcommand{\Ac}{\mathcal{A}}

\newcommand{\Cc}{\mathcal{C}}

\newcommand{\Lc}{\mathcal{L}}

\newcommand{\Mc}{\mathcal{M}}

\newcommand{\Tc}{\mathcal{T}}

\newcommand{\CCC}{\mathfrak{C}}

\newcommand{\LLL}{\mathfrak{L}}
\newcommand{\GGG}{\mathfrak{G}}
\newcommand{\III}{\mathfrak{I}}
\newcommand{\SSS}{\mathfrak{S}}
\newcommand{\SSSt}{\mathfrak{St}}

\title{Stacks on stratified spaces}

\begin{document}
\maketitle
\begin{abstract}
In this paper, we go into the study of the $2$-category $\SSS_{\Sigma}$ of $\Sigma$-cons\-truc\-tible stacks . We show the $2$-equivalence between $\SSS_{\Sigma}$ and a combinatoric $2$-category whose objects are given by a $2$-representation of each stratum plus some gluing data.

\end{abstract}

A {\bf stack} is a generalization of a sheaf of categories. The notion of equality between two categories being too strict, a stack is, roughly speaking, a ``sheaf of categories up to equivalence''. This lax version of sheaves allows to give a structure to objects that can be glued up to isomorphisms. For example, Beilinson, Bernstein and Deligne have shown in \cite{BBD} that if $X$ is a topological space, the data for all open $U$ of $X$, of the category of perverse sheaves on $U$ is a stack.

Most of the notations and properties of sheaves can be extended to
stacks. In this article we focus on the notions of locally constant
stack and constructible stack on a stratified space.

Let us recall some classical facts on sheaf theory. It is well knowns
that the category of locally constant sheaves on a locally 1-connected
topological space $X$ is equivalent to the category of representations
of the fundamental groupoid $\Pi_1(X)$. This result gives a
topological description of the category of locally constant sheaves on
$X$. Constructible sheaves are a natural generalization of locally
constant sheaves. Indeed, a sheaf $\mathcal F$ is constructible if
there exists a stratification $\{\Sigma_k\}$ of $X$ such that, for every
stratum $\Sigma_k$, the restriction of $\mathcal F$ to $\Sigma_k$ is
locally constant on $\Sigma_k$. Another classical construction for
sheaves on a topological space $X$, is the description of a sheaf
through some descent data. This gives an answer to
the natural question of how we can recover a sheaf $\mathcal F$ from its restrictions to
open or closed subset of $X$. In particular, $\mathcal F$ is uniquely
determined by its restrictions to an open set $U\subset X$ and its
complementary $F$, plus some gluing data given by the morphisms of
adjunction
$$ i^{-1}_F\mathcal F\to i^{-1}_Fi_{U*}i^{-1}_U \mathcal F \ ,$$
where $i_F$ and $i_U$ are the inclusions of $F$ and $U$ in
$X$. Combining the topological description of locally constant sheaves
given above with the previous gluing construction, one can obtain a
combinatorial description of a constructible sheaf with respect to a
stratification.

In \cite{Ingo2}, P. Polesello and I. Waschkies generalized to the $2$-category
of locally constant staks on a suitable topological space the
topological description of locally constant sheaves cited above. In
particular they introduced the $2$-monodromy functor from the $2$-category of locally constant stacks to the
$2$-category of $2$-representations of $\Pi_2(X)$. In the \textbf{first
  section} of this paper we recall the definition of the $2$-monodromy
functor and defining a quasi-$2$-inverse slightly different from the one
in \cite{Ingo2}. Then, given a locally trivial fiber bundle $p:X\to B$, we
consider the $2$-functors of direct and inverse image relative to $p$
between the $2$-categories of locally constant stacks and we translate
such functors in the language of $2$-representations. 

In the \textbf{second section}, we consider a stratified topological
space and we study how to recover a stack from its restrictions to the
strata. In particular, we generalize the gluing construction cited above for sheaves
to the case of stacks. Given a topological space $X$ and a
stratification $\Sigma=\{\Sigma_k\}_k$ of $X$, we define the
$2$-category $\mathfrak S_\Sigma$ whose objects are given by
\begin{itemize}
\item a stack $\CCC_k$ on each stratum $\Sigma_k\underset{i_k}\hookrightarrow X$, 
\item a functor of stacks 
$F_{kl} : \CCC_{k}\longrightarrow i_{k}^{-1}i_{l*}\CCC_{l}$,
for every couple $\Sigma_{k}, \Sigma_{l}$ of strata such that $\Sigma_{k}\subset \overline{\Sigma}_{l}$,
\item for every triple $\Sigma_k, \Sigma_l, \Sigma_m$ such that $\Sigma_{k}\subset \overline{\Sigma}_{l}\subset \overline{\Sigma}_m$, some morphisms of functors. 
\end{itemize}

We show the following 

\begin{thm}
  The $2$-category $\mathfrak St_X$ of stacks on $X$ is equivalent to
  the $2$-category $\mathfrak S_\Sigma$.
\end{thm}

Hence we see that, in order to define a stack on a stratified
topological space, it is sufficient to have stacks on each stratum plus some gluing data consisting of functors of stacks and morphisms of functors. To prove the theorem, we define a couple of quasi-$2$-inverse functors : the ``restriction functor'' $R_{\Sigma}$ going from $\SSS t_{X}$ to $\SSS_{\Sigma}$ and the ``gluing functor'' $G_{\Sigma}$. The former is the restriction of a stack to each stratum, plus some functors and morphisms given by the $2$-adjunction between $i_{k*}$ and $i_{k}^{-1}$. The definition of the latter is more technical. For all object of $\SSS_{\Sigma}$ we define a $2$-functorial $2$-limit encoding the gluing data. \\

In the third section we focus on \textbf{constructible stacks}. The
notions of constructible stack was introduced by D. Treumann in \cite{Tr1}. It is a natural generalization of constructible sheaf. A stack $\CCC$ is called constructible if there exists a stratification $\Sigma$ of $X$ such that $\CCC$ is locally constant along each stratum.  In \cite{Tr1}, D. Treumann has also introduced the exit-path $2$-category, which is a stratified version of the fundamental 2-groupoid and  he showed that these two $2$-categories are equivalent. Let us also cite J. Woolf in \cite{Wo}, he generalize the work of  D. Treumann to homotopically stratified sets.

In what follows we focus on constructible stacks with respect to a
fixed stratification $\Sigma$ of $X$. Although we are
interested in the same $2$-category $\SSS t_{\Sigma}^{c}$ of $\Sigma$-constructible stacks, our approach is different. We show the $2$-equivalence between $\SSS t_{\Sigma}^{c}$ and a $2$-category whose objects are combinatoric data of $2$-re\-pre\-sen\-ta\-tions, functors of $2$-representations and isomorphisms of functors. As a constructible stack is locally constant along each stratum  and as the $2$-monodromy defined by P. Polesello and I. Waschkies is an equivalence of categories, it is natural to ask if the data for every stratum $\Sigma_{k}$ of a $2$-representation of $\Pi_{2}(\Sigma_{k})$ is sufficient to define a unque constructible stack, up to equivalence.

Now, if we want to describe combinatorially the $2$-category $\SSS t_{\Sigma}^c$ of $\Sigma$-cons\-truc\-tible stacks, it remains to understand how the gluing data can be read in the language  of $2$-representations. To have a better understanding we restrict ourself to the case of Thom-Mather spaces. A Thom-Mather
space is a stratified space plus a tubular neighborhood $T_k$ of each stratum $\Sigma_k$ together with a locally trivial fiber bundle $p_{k} : T_{k}\rightarrow \Sigma_{k}$, (for precise definitions see \cite{Ma} and \cite{T}). In this case, we show that the $2$-functor $i_{k}^{-1}i_{l*}$ restricted to the $2$-category of locally constant stacks on $\Sigma_{l}$ is $2$-equivalent to the functor $p_{k*}i_{kl}^{-1}$,  where $i_{kl}$ is the  natural inclusion of $\Sigma_{l}\cap T_{k} $ in $\Sigma_{l}$. Now, in the first section, we have defined the equivalent functor in the $2$-category of $2$-representations. Hence we define a $2$-category $\SSS t_{\Sigma}^c$, $2$-equivalent to the $2$-category of constructible stacks, 
whose objects are given by:
\begin{itemize}
\item for every stratum $\Sigma_{k}$, a $2$-representation $\alpha_{k}$ of the fundamental $2$-groupoid $\Pi_{2}(\Sigma_{k})$,
\item for every couple $\Sigma_{k}$ and $\Sigma_{l}$ of strata such that $\Sigma_{k}\subset \overline{\Sigma}_{l}$, a functor of $2$-representation $F_{kl}$ : 
$$F_{kl} : \alpha_{k}\longrightarrow p_{k*}i_{kl}^{-1}\alpha_{l} $$
\item some morphisms of functors.
\end{itemize}

As a $2$-representation of a $2$-groupoid is equivalent to the data of categories, functors of categories and isomorphisms of functors, we can conclude our combinatorial description of a constructible stack.

As an application , in \cite{Tr2},
D. Treumann has used his description of the $2$-category of constructible stacks and a description of the category of perverse sheaves given by MacPherson and Vilonen in \cite{McV} to characterize the stack of perverse sheaves and in the case of Thom-Mather spaces he has showed that if the stratum are $2$-connected the category of perverse sheaves is equivalent to the category of finite-dimensional modules over a finite-dimensional algebra. As he has used a
non explicit local description he does not obtain an explicit description. In the same spirit, using the description of the category $\mathfrak S_\Sigma$, we glue I glue in  \cite{these}   descriptions of the category of perverse sheaves on a normal crossing given by A. Galligo, M. Granger and Ph. Maisonobe in \cite{GGM}, to obtain  explicit descriptions of the category of perverse sheaves on smooth toric varieties stratified by the torus action. For a presentation of the result see \cite{vari}.
~\\

{\bf Conventions}. Here we use the term ``$2$-category'' for a strict $2$-category. It means that the composition of $1$-morphisms is strictly associative. By a $2$-functor, we mean a morphism of $2$-category preserving the composition of $1$-morphisms up to isomorphism. By a $2$-representation of a $2$-groupoid $G$, we mean a $2$-functor from $G$ to the $2$-category of categories $\Cc\Ac\Tc$. If $\alpha$ is a $2$-representation of athe fundamental groupoid $\Pi_2(X)$ of $X$ and $F$ is a subset of $X$, with an abusive notation, we also denote $\alpha$ the functor $\alpha$ restricted to $F$.

We do not recall the definitions of stack, constant stack and locally constant stack, the reader can find them in \cite{Ingo2} or in \cite{Tr1}. We often use the notion of $2$-limit and $2$-colimit, their definition is given in \cite{SGA4-2}, for an explicit description see  for example the annex of \cite{Ingo1} or \cite{Interchange}. The $2$-adjunction plays an important role in this paper, we refer to \cite{Gray}.

\section{Locally constant and constructible stacks}
Let $X$ be a locally connected space. \\
In this section we go into the study of the equivalence between the $2$-category $\LLL_{X}$ of locally constant stacks on $X$ and the $2$-category, $Rep(\Pi_{2}(X), \Cc\Ac\Tc)$ of $2$-representations of the fundamental $2$-groupoid $\Pi_{2}(X)$ of $X$. 

In a first time we shortly recall the definition of the $2$-monodromy given by  P. Polesello and I. Waschkies in \cite{Ingo2}. This $2$-equivalence, denoted $\mu$, is a generalization of the monodromy going from the category of locally constant sheaves on $X$ to the category of representation of the fundamental groupoid of $X$. They show that $\mu$ is an equivalence defining a quasi-$2$-inverse. Here we define a quasi-$2$-inverse $\nu$ of $\mu$ slightly different from the one given in \cite{Ingo2}. 

Then, we translate in the language of $2$-representations some operations on locally constant stacks. More precisely, let $f: Y\rightarrow X$ be a continuous map and $p: X\rightarrow B$ be a locally trivial fiber bundle. If $\CCC$ is a locally constant stack, then $f^{-1}(\CCC)$ and $p_{*}\CCC$ are locally constant. 
We define  two $2$-functors, also denoted $f^{-1}$ and $p_{*}$,  going from $2Rep(\Pi_{2}(X), \Cc \Ac\Tc)$ to  $2Rep(\Pi_{2}(Y), \Cc \Ac\Tc)$ and $2Rep(\Pi_{2}(B), \Cc\Ac\Tc)$ respectively, commuting with the $2$-monodromy.

~\\\\
Let $\CCC$ be a locally constant stack on $X$. Let $\gamma : I \rightarrow X$ be a path in $X$.\\ 
As $I$ is contractible, the stack $\gamma^{-1}\CCC$ is a constant stack, thus the following functors are equivalences :
$$\CCC_{x_{0}} \simeq (\gamma^{-1}(\CCC))_{0} \buildrel\sim\over\leftarrow \Gamma([0,1], \gamma^{-1}(\CCC))\buildrel\sim\over\rightarrow (\gamma^{-1}(\CCC)_{1}\simeq \CCC_{x_{1}}.$$
Let us denote $\tilde{\gamma}$ the composition of the previous equivalences. \\
Let $H$ be a homotopy in $X$, going from a path $\gamma_{0}$ to a path $\gamma_{1}$.
The following diagram commutes up to isomorphisms. 
$$\xymatrix@!0 @R=0.7cm @C=1.2cm{
\CCC_{x_{0}}\eq[dr] \ar@/^{0.7cm}/[rrrrrrrr]^{\tilde{\gamma}_{0}} \ar@/_{0.7cm}/[dddddddd]_{Id}&&&&&&&&\CCC_{x_{1}} \eq[dl] \ar@/^{0.7cm}/[dddddddd]^{Id}\\ 
&(H^{-1}\CCC)_{(0,0)}  &&& \ar[lll]_\sim\Gamma\big(\{0\}\times I, H^{-1}\CCC\big) \ar[rrr]^\sim &&& (H^{-1}\CCC)_{(0,1)} \\\\\\
&\Gamma\big(I\times \{0\}, H^{-1}\CCC\big) \ar[uuu]^\sim \ar[ddd]_{\sim} &&& \Gamma\big(I\times I, H^{-1}\CCC\big) \ar[rrr]^\sim  \ar[uuu]^\sim \ar[ddd]_{\sim} \ar[lll]_{\sim} &&& \Gamma\big(I\times \{1\}, H^{-1}\CCC\big) \ar[uuu]_{\sim} \ar[ddd]^\sim \\\\\\
&(H^{-1}\CCC)_{(1,0)}  &&&\Gamma\big(\{1\}\times I, H^{-1}\CCC\big) \ar[rrr]_{\sim} \ar[lll]^\sim &&& (H^{-1}\CCC)_{(1,1)}\\
\CCC_{x_{0}} \eq[ur] \ar@/_{0.7cm}/[rrrrrrrr]_{\tilde{\gamma}_{1}}&&&&&&&&\CCC_{x_{1}} \eq[ul]
}$$
The suitable composition of previous isomorphisms of functors gives an isomorphism of functors :
$$\tilde{H} : \tilde{\gamma}_{0} \longrightarrow \tilde{\gamma}_{1}.$$
Then, the image of $\CCC$ by the $2$-monodromy $\mu$ is the $2$-functor defined as follow : 
$$\begin{array}{cccc}
\Pi_{2}(Y) & \longrightarrow & \Cc \Ac\Tc \\
x & \longmapsto & \CCC_{x}\\
\gamma : x_{0}\rightarrow x_{1} & \longmapsto & \tilde{\gamma} :\CCC_{x_{0}} \buildrel\sim\over\rightarrow \CCC_{x_{1}}\\
H : \gamma_{0}\rightarrow \gamma_{1} & \longmapsto & \tilde{H} : \tilde{\gamma}_{0}\rightarrow \tilde{\gamma}_{1}
\end{array}$$
As the equivalences, the isomorphisms of functors come from a $2$-functor this application is $2$-functorial.\\
Now, let us define a quasi-$2$-inverse of the $2$-monodromy, denoted $\nu$. Let $\alpha$ be a $2$-representation of $\Pi_{2}(X)$, let us consider :
\begin{itemize}
\item for all open $U$ of $X$, the category $\nu(\alpha)(U)=\begin{displaystyle}2\!\!\!\varprojlim_{ \Pi_{2}(U)}\alpha\end{displaystyle}$, 
\item for every pair $V\subset U$ of open subsets of $X$, the functor :
$$\nu(\alpha)(U)\begin{displaystyle} 2\!\!\!\varprojlim_{ \Pi_{2}( U)}\alpha \longrightarrow 2\!\!\!\varprojlim_{ \Pi_{2}( V)}\alpha  \end{displaystyle}=\nu(\alpha)(V)$$
defined by the projections $\pi_{x} : \begin{displaystyle} 2\!\!\!\varprojlim_{ \Pi_{2}( U)}\alpha  \end{displaystyle}\rightarrow \alpha(x)$ of the $2$-limit.
\item for every triple $W\subset V\subset U$, the isomorphism of functors defined by the isomorphisms given by the $2$-limit :
$$\xymatrix{\begin{displaystyle} 2\varprojlim_{ \Pi_{2}(U)}\alpha \end{displaystyle} \ar[r] \ar[rd] & \begin{displaystyle} 2\varprojlim_{ \Pi_{2}(V)}\alpha\end{displaystyle} \ar[d]\\
& \begin{displaystyle} 2\varprojlim_{ \Pi_{2}(W)}\alpha\end{displaystyle}}$$
\end{itemize}
\begin{lem}
These data define a locally constant stack in a $2$-functorial way.
\end{lem}
\dem
The proof is similar to the proof of the theorem 2.2.5 of \cite{Ingo2}.
\cqfd

\begin{thm}
Let $X$ be a relatively $2$-connected space, $\Pi_{2}(X)$ its fundamental $2$-groupoid. Then the $2$-functors $\mu$ and $\nu$ are $2$-equivalent.
\end{thm}
It would be interesting to establish a dictionary between the operations on locally constant stacks and operations on the $2$-representations of $\Pi_{2}$. \\
The inverse image of a locally constant stack by a continuous function is locally constant.
\begin{prop}\label{restriction}
 Let $f: X\rightarrow Y$ be a continuous map. Let $f^{-1}$ denote the $2$-functor defined by  :
$$\begin{array}{cccc}
f^{-1} :& Rep\big(\Pi_{2}(Y), \Cc\Ac\Tc\big) &\longrightarrow & Rep\big(\Pi_{2}(X), \Cc\Ac\Tc\big)\\
& \alpha & \longmapsto& \left(\begin{array}{ccccccc}
                                              x\in X & \mapsto & \alpha(f(x))\\
                                               \gamma : I\rightarrow  X & \mapsto & \alpha(f\circ \gamma)\\
                                               \varepsilon: I\times I\rightarrow  X & \mapsto & \alpha(f\circ \varepsilon)
                                               \end{array}\right)
\end{array}$$
Then, the two $2$-functors $f^{-1}$ and $\mu\circ f^{-1}\circ \nu$ going from $Rep\big(\Pi_{2}(Y), \Cc\Ac\Tc\big)$ to $Rep\big(\Pi_{2}(X), \Cc\Ac\Tc\big)$,  are equivalent.
\end{prop}
\dem
The proof is straightforward.
\cqfd
Now, we are interested in the direct image of a locally constant stack. Not all direct images of locally constant stack is a locally constant stack, but D. Treumann has shown the following proposition : 
\begin{prop}[\cite{Tr1}]
Let $X$ and $B$ be locally contractible spaces. Let $p :X  \rightarrow B$  be a locally
trivial fiber bundle. Let $\CCC$ be a locally constant stack on $X$, then $p_{*}\CCC$ is locally constant on $B$ .
\end{prop}
Let $X$ and $B$ be two locally contractible spaces and $p$ be a locally trivial fiber bundle : 
$$p: X\rightarrow B.$$  
In what follows, we define explicitly a functor $p_{*}$ going from \linebreak $2Rep(\Pi_{2}(X), \Cc \Ac\Tc)$ to $2Rep(\Pi_{2}(B), \Cc at)$ such that the following diagram commutes up to isomorphisms  :
$$\xymatrix @!0 @C=3.2pc @R=2pc{\cat L_{X} \ar[rrr]^{p_{*}}  \ar[ddd]_{\mu} &&&\cat L_{B} \ar[ddd]^{\mu} \\
 && \ar@{=>}[ld]_{\sim}\\
& & &\\
2Rep(\Pi_{2}(X), \Cc at) \ar[rrr]_{p_{*}} &&& 2Rep(\Pi_{2}(B), \Cc at).
}$$
But first we need to fix some data and notations. \\
In all this section $\alpha$ denotes a $2$-representation of the fundamental $2$-groupoid $\Pi_{2}(X)$. 
If $x\in Ob(\Pi_{2}(X))$ we denote by $\alpha_{x}$ the isomorphism of functors : 
$$\alpha_{x} : \alpha(Id_{x})\buildrel\sim\over\longrightarrow Id_{\alpha(x)}.$$
If $\gamma$ and $\gamma'$ are two composable paths in $X$, we denote by $\alpha_{\gamma, \gamma'}$ the isomorphism~:
$$\alpha_{\gamma,\gamma'}: \alpha(\gamma)\circ\alpha(\gamma')\buildrel\sim\over\longrightarrow \alpha(\gamma\land\gamma').$$
Let  $\gamma:b_{0}\rightarrow b_{1}$ be a path in $B$ and $\varepsilon: \gamma_{0}\rightarrow \gamma_{1}$ be a homotopy in $B$, for $t\in I$ and for  $t_{1}, t_{2}\in I^2$ we set :
$$F_{t}:=p^{-1}\big(\gamma(t)\big),~~~~F_{(t_{1},t_{2})}:=p^{-1}\big(\varepsilon(t_{1}, t_{2})\big).$$
Because of the contractibility of $I$ and $I\times I$, the fibrations $p_{1\gamma}$ and $p_{1\varepsilon}$ given by the pullbacks 
$$
\xymatrix{
 I\times_{\gamma}X\ar[r]^{p_{2\gamma}}\ar[d]_{p_{1\gamma}} &X\ar[d]^{p}\\
 I\ar[r]_{\gamma} &B
}~~~~
\xymatrix{
 (I\times I)\times_{\varepsilon}X\ar[r]^{~~~~~~~~p_{2\varepsilon}}\ar[d]_{p_{1\varepsilon}} &X\ar[d]^{p}\\
 I\times I\ar[r]_{\varepsilon} &B
}
$$
are trivialisable, and there exist trivializations $h_{\gamma}$ and $h_{\varepsilon}$ of $I\times_{\gamma}X$ and \linebreak $(I\times I)\times_{\varepsilon}X$  
$$
\xymatrix{
I\times F_{0}\ar[r]^{h_{\gamma}}\ar[rd]_{\pi} & I\times_{\gamma}X\ar[r]^{p_{2}}\ar[d]^{p_{1}} &X\ar[d]^{p}\\
& I\ar[r]_{\gamma} &B
}~~~~~~~
\xymatrix{
(I\times I)\times F_{0}\ar[r]^{h_{\varepsilon}}\ar[rd]_{\pi} & I\times_{\varepsilon}X\ar[r]^{p_{2\varepsilon}}\ar[d]^{p_{1}} &X\ar[d]^{p}\\
& I\ar[r]_{\varepsilon} &B
}
$$
such that  
\begin{equation}\label{equation_projection}p_{2}\circ h_{\gamma}\mid_{\{0\}\times F_{0}}=\pi_{2} \text{~~and~~}p_{2}\circ h_{\varepsilon}\mid_{\{0\}\times\{0\}\times F_{0}}=\pi_{3}.\end{equation}
These trivializations are not unique but two such trivializations are homotopic.
\begin{lem}\label{trivialization}
If $h$ and $h'$ are two trivializations of $I\times_{\gamma}X$ satisfying the condition (\ref{equation_projection}) then $h_{\gamma}$ and $h_{\gamma}'$ are homotopic and there exists a  homotopy H, unique up to homotopy, between them such that :
\begin{equation}\label{projection2}(p_{2}\circ H)\mid_{I\times\{0\}\times X}=\pi_{3},~~~~~~~~p_{1}\circ H=\pi_{2},\end{equation}
and such that, for all $t\in I$ $H(t, \cdot, \cdot)$ is an isomorphism  from $I\times F_{0}$ to $I\times_{\gamma} X$.\\
In the same way, if $h_{\varepsilon}$ and $h_{\varepsilon}'$ are two trivializations of $(I\times I)\times_{\varepsilon}F_{0}$, then they are homotopic. 
\end{lem}
\dem
We set  $h^{-1}=\big((h^{-1})_{1},(h^{-1})_{2}\big)$. \\
Let us remark that, with the above notations $(h_{1}^{-1})_{1}=p_{1}$ and as $h$ and $h'$ satisfy the condition \ref{equation_projection}, for all $x\in F$, we have : $$h(0,x)=h'(0,x).$$
This assures that the application 
$$\begin{array}{cccc}
H : & I\times I\times F_{0} & \longrightarrow & I\times_{\gamma}X\\
& (t_{1},t_{2},x)& \longmapsto & h\Big(t_{2},(h^{-1})_{2}\big(h'(t_{1}t_{2},x)\big)\Big)
\end{array}$$
is a homotopy from $h$ to $h'$ and moreover that the conditions of the lemma are satisfied. \\
Now, let us suppose that $H_{1}$ and $H_{2}$ are two such homotopies going from 
$h$ to $h'$, satisfying the conditions (\ref{projection2}). Let $t\in I$, let us set :
$$H_{1}(t,\cdot, \cdot)^{-1}:=\Big(H_{1}(t,\cdot, \cdot)^{-1}_{1}, H_{1}(t,\cdot, \cdot)^{-1}_{2}\Big).$$
As above, the condition (\ref{projection2}) assures that for all $(t,x)\in I\times F_{0}$ we have the equality 
$$H_{1}(t,0,x)=H_{2}(t,0,x).$$
Hence the  homotopy 
$$\begin{array}{cccc}
I\times I\times I\times F_{0} &\longrightarrow & I\times_{\gamma}X\\
(t_{1}, t_{2}, t_{3}, x)&\longmapsto & H_{1}\Big(t_{2},t_{3}, \big(H_{1}(t_{2}, \cdot ,\cdot) )^{-1}_{2}\circ H_{2}(t_{2}, t_{1}t_{3}, x)\Big)
\end{array}$$
goes from $H_{1}$ to $H_{2}$.\\
If $\varepsilon$ is a homotopy in $B$, we define in the same way the homotopy between two trivializations of $(I\times I)\times_{\varepsilon}X$.
\cqfd
These trivializations and their unicity up to homotopy allow us to use the properties of the product space, in particular the fact that $$\pi_{n}(X\times Y)\simeq \pi_{n}(X)\times \pi_{n}(Y).$$ 
Let $\gamma :x_{0}\rightarrow x_{1}$ be a path in $B$, we set $F_{0} :=p^{-1}(\gamma(0))$ and $F_{1}=p^{-1}(\gamma(1))$.
In what follows, to each trivialization of $I\times_{\gamma}X$ satisfying the condition (\ref{equation_projection}), we associate a functor 
 $$\begin{displaystyle}2\!\!\!\! \varprojlim_{\Pi_2(F_{0})}\alpha \end{displaystyle}\longrightarrow \begin{displaystyle}2\!\!\!\! \varprojlim_{\Pi_2(F_{1})}\alpha \end{displaystyle},$$and  to each homotopy between two such trivializations we define a morphism of functors between the functors : 
$$\xymatrix{\begin{displaystyle} 2\varprojlim_{\Pi_{2}(F_{0})}\end{displaystyle} \UN[rr]{}{}{} && \begin{displaystyle}2\varprojlim_{\Pi_{2}(F_{1})} \end{displaystyle}}.$$
Let $h$ be a trivialization of $I\times_{\gamma}X$, we denote by $\Gamma_{h}$  the composition $\Gamma_{h}=p_{2}\circ h.$
Hence the following diagram commutes  
\begin{equation}\label{diag}
\xymatrix{
I\times F_{0} \ar[r]^{\Gamma_{h}} \ar[d]_{p_{1}} & X\ar[d]^{p} \\
I\ar[r]_{\gamma} & B. 
}\end{equation}
Let us remark that for all  $t\in I$,  the  application 
$$\begin{array}{cccccccccc}
\Gamma_{h}(t, \cdot) : &F_{0}& \longrightarrow &F_{t}\\
& x&\longmapsto & \Gamma(t,x)
\end{array}$$
is an isomorphism and for all $x\in F_{0}$, the application 
$$\begin{array}{cccccccccc}
\Gamma_{h}(\cdot, x) : &I& \longrightarrow &X\\
& t&\longmapsto & \Gamma_{h}(t,x)
\end{array}$$
is a path from $x$ to $\Gamma(1, x)$.
\begin{lem}\label{functor}
Let $h$ be a trivialization of $I\times_{\gamma} X$. With the above notations, the data  of:
\begin{itemize}
\item for every $y\in F_{1}$,  the functors $\alpha\big(\Gamma_{h}(\cdot,x)\big)\circ \pi_{x}$, where $y=\Gamma_{h}(1,x)$,
\item for every path $\delta_{1} : x_{1}\rightarrow y_{1}$ in $F_{1}$,  the morphism of functors \linebreak $\big(Id\bullet \pi_{\delta_{0}}\big)\circ\big(\alpha(\Gamma_{h}\circ H_{\delta_{1}})\bullet Id\big)$ visualized by :
$$\xymatrix @!0 @C=3.2pc @R=2pc{
&&&&\alpha(x_{0}) \ar[rrr]^{\alpha\big(\Gamma_{h}(\cdot,x_{0})\big)}  \ar[dddd]_{\alpha\big(\delta_{0}\big)} &&&\alpha(x_{1}) \ar[dddd]^{\alpha(\delta_{1})} \\
&&&\ar@{=>}[ldd]_{\sim}^{\pi_{\delta_{0}}}& && \ar@{=>}[ldd]_{\sim}^{\alpha\big(\Gamma_{h}(H_{\delta_{0}})\big)}\\
&\begin{displaystyle}2\!\!\!\varprojlim_{\Pi_{2}(F_{0})}\alpha \end{displaystyle} \ar@/_{0.5cm}/[rrrdd]_{\pi_{y}} \ar@/^{0.5cm}/[rrruu]^{\pi_{x}}&&& &&&&\\
&&&&& & &\\
&&&& \alpha\big(y_{0}\big) \ar[rrr]_{\alpha\big(\Gamma_{h}(\cdot,y_{0})\big)} &&& \alpha\big(y_{1}\big)
}$$ 
 where :
\begin{itemize}
\item  $\pi_{x}$, $\pi_{y}$ and $\pi_{\delta_{0}}$  are  the functors and the isomorphisms of functors given by the $2$-limit,
\item the path $\delta_{0}$, is the inverse image of $\delta_{1}$, by $\Gamma_{h}(1, \cdot)$, i.e. $\delta_{1}=\Gamma_{h}(1,\delta_{0}),$
\item and $H_{\delta_{0}}$ is a homotopy in $I\times F_{0}$ between the paths  $(1,\delta_{0})\land(Id,x_{0})$ and  $(Id,y_{0})\land(0,\delta_{0})$. 
\end{itemize} 
 \end{itemize}
 define a functor, also denoted $\Gamma_{h}$, 
$$\begin{displaystyle}\Gamma_{h} : 2\!\!\!\!\varprojlim_{\Pi_{2}( F_{0})}\alpha\longrightarrow 2\!\!\!\!\varprojlim_{\Pi_{2}( F_{1})}\alpha\end{displaystyle}.$$
\end{lem}
\dem We have to show that these data satisfy the commutation conditions. \\
Let $\delta_{1} : x_{1}\rightarrow y_{1}$ and $\delta'_{1} : y_{1}\rightarrow z_{1}$ two composable paths in $F_{1}$. We denote by $x_{0}$, $y_{0}$ and $z_{0}$ the points of $F_{0}$ such that 
$$x_{1}=\Gamma(1, x_{0}), ~~~y_{1}=\Gamma(1, y_{0}), ~~~z_{1}=\Gamma(1, z_{0})$$
and $\delta_{0}$, $\delta_{0}'$ the paths in $F_{0}$ such that : 
$$\delta_{1}=\Gamma(1, \delta_{0}), ~~~\delta_{1}'=\Gamma(1, \delta_{0}').$$
In view of the commutation conditions satisfied by the functors $\pi_{x}$ and the morphisms of functors $\pi_{\delta}$, we have to show that :
$$\alpha(\Gamma_{h}\circ H_{\delta_{1}\land\delta'_{1}})=\alpha\Big(((\Gamma_{h}\circ H_{\delta_{1}'})\bullet Id)\circ(Id\bullet(\Gamma_{h}\circ H_{\delta_{1}}))\Big).$$
As $\alpha$ in a $2$-representation of the $2$-groupoid $\Pi_{2}(X)$ it is sufficient to show that there exists a homotopy from $\Gamma_{h}\circ H_{\delta_{1}\land\delta_{1}'}$ to $((\Gamma_{h}\circ H_{\delta_{1}'})\bullet Id)\circ(Id\bullet(\Gamma_{h}\circ H_{\delta_{1}}))$. Let us first consider the two homotopies in $I\times F_{0}$,  $H_{\delta_{1}\land\delta_{1}'}$ and $(H_{\delta_{1}'}\bullet Id)\circ (Id\bullet H_{\delta_{1}})$, as $\pi_2(I\times F_0) \simeq \pi_2(I)\times \pi_2(F_0)$, there exists an homotopy between them. Applying $\Gamma_{h}$ we find the homotopy we looked for. \\
We use the same arguments to show that if $\varepsilon : \delta\rightarrow \delta'$ is a homotopy in $F_{1}$  we have the equality : 
$$\alpha(\varepsilon)\circ (Id\bullet \pi_{\delta})\circ (\alpha(\Gamma_{h}\circ H_{\delta})\bullet Id)=(Id\bullet \pi_{\delta'})\circ (\alpha(\Gamma_{h}\circ H_{\delta'})\bullet Id).$$
Hence the data given in the statement of the lemma define a functor from $\begin{displaystyle} 2\!\!\!\!\varprojlim_{\Pi_{2}( F_{0})}\alpha\end{displaystyle}$ to $\begin{displaystyle} 2\!\!\!\!\varprojlim_{\Pi_{2}( F_{1})}\alpha\end{displaystyle}$.
\cqfd
Now let us consider two trivializations $h_{0}$ and $h_{1}$ of $I\times_{\gamma}X$ satisfying the conditions (\ref{equation_projection}), and $H: h_{0}\rightarrow h_{1}$  a homotopy between them satisfying the conditions given in the lemma \ref{trivialization}. We define an isomorphism of functors from $\Gamma_{h_{0}}$ and $\Gamma_{h_{1}}$. In what follows, we set :
$$\Gamma_{0}:=\Gamma_{h_{0}}~~~~\text{and}~~~~\Gamma_{1}:=\Gamma_{h_{1}}.$$
Let $x_{1}\in F_{1}$. We denote by $x_{0}$ and $y_{0}$ the points such that 
$$\Gamma_{0}(1,x_{0})=x_{1}~~~\text{and}~~~\Gamma_{1}(1,y_{0})=x_{1}.$$
Let us remark that the application 
$$\begin{array}{cccc}
\beta_{x_{1}} :& I&\longrightarrow &F_{0}\\
& t&\longmapsto & (\Gamma_{1}(1, \cdot))^{-1}\circ H(t,1,x_{0})
\end{array}$$
is a path in $F_{0}$ going from $y_{0}$ to $x_{0}$.
Let us consider a map 
$$\begin{array}{cccccc}
\phi : & I\times I & \longrightarrow & I\times I\\
&  (t_{1}, t_{2}) & \longmapsto & \left\{ \begin{array}{cccc}
                                                 (2t_{2}, 0) & \text{if}&  0\leq t_{2}\leq \dfrac{t_{1}}{2}\\
                                                 (t_{1}, \dfrac{2}{2-t_{1}}t_{2}-\dfrac{t_{1}}{2}) &\text{if}&\dfrac{t_{1}}{2}\leq t_{2}\leq 1.
                                                 \end{array} \right.
\end{array}$$
Then the composition, denoted $H_{x_{1}}$, of $\phi$  with the following map:
$$\begin{array}{cccc}
 I\times I& \longrightarrow & X\\
  (t_{1},t_{2}) &\longmapsto & ((\Gamma_{1}(t_{2},\cdot)^{-1}\circ H(t_{1}, t_{2},x_{0})
\end{array}$$
is a homotopy from $\Gamma_{0}(\cdot , x_{0})\land \beta_{x_{1}}$ to $\Gamma_{1}(\cdot, y_{0})$. 
\begin{lem}\label{isomorphism}
With the above notation, the data for every $x_{1}\in F_{1}$ of the suitable composition of isomorphisms visualized by the following diagram : 
$$\xymatrix @!0 @C=1.7pc @R=1.3pc{
&&&&&\alpha(x_{0}) \ar@/^{0.5cm}/[rrrrrddd]^{\alpha(\Gamma_{0}(\cdot, x_{0}))} \ar[dddddd]_{\alpha(\beta_{x_{1}})}&&& \\
&&&&&&&&\\
&&&&\ar@{=>}[lldd]^{\sim}_{\pi_{\beta_{x_{1}}}}& &&&\ar@{=>}[lldd]^{\sim}_{\alpha(H_{x_{1}})}&&& \\
\begin{displaystyle}2\!\!\!\varprojlim_{\Pi_{2}(F_{0})}\alpha \end{displaystyle} \ar@/_{0.5cm}/[rrrrrddd]_{\pi_{y_{0}}} \ar@/^{0.5cm}/[rrrrruuu]^{\pi_{x_{0}}}&&&&& && &&&  \alpha(x_{1})  \\
&&&&&&&&&&\\
&&&&& &&&&&\\
&&&&& \alpha(y_{0}) \ar@/_{0.5cm}/[rrrrruuu]_{\alpha(\Gamma_{1}(\cdot, y_{0}))}
}$$ 
define an isomorphism of functors also denoted $H$ :
$$\xymatrix{
\begin{displaystyle}2\!\!\!\varprojlim_{\Pi_{2}(F_{0})}\alpha \end{displaystyle}\UN[rrr]{\Gamma_{0}}{\Gamma_{1}}{H}&&&\begin{displaystyle}2\!\!\!\varprojlim_{\Pi_{2}(F_{1})}\alpha \end{displaystyle}
}$$
\end{lem}
\dem 
The proof uses the same arguments as the proof of lemma \ref{functor}. Hence if $\delta$ is a path in $F_{1}$, using the facts that $\pi_\beta$ satisfied commutation conditions, there exists of homotopy of homotopies
$$\begin{array}{cccc}
 I\times I \times I :& \longrightarrow & X\\
(t_{1}, t_{2}, t_{3}) &\longmapsto & H\big(t_{1}, H_{\delta}(t_{2}, t_{3})\big)
\end{array}$$ 
and the uniqueness up to homotopy of $H_{\delta}$,
we have that the data given in the lemma satisfy the commutations conditions and define an isomorphism from $\Gamma_{0}$ to $\Gamma_{1}$.
\cqfd
Let $\gamma_{0}, \gamma_1:x_{0}\rightarrow x_{1}$ a be two paths in $B$, $\varepsilon : \gamma_{0}\longrightarrow \gamma_{1}$ a homotopy in $B$ and $h_{\varepsilon}$ a trivialization of $(I\times I)\times_{\varepsilon}X$. 
As above, we denote by $\Gamma_{\varepsilon}$ the composition $p_{2}\circ h_{\varepsilon}$. For $i=0$ or $i=1$, we set : 
$$\Gamma_{\varepsilon i}(t,x):=\Gamma_{\varepsilon}(i,t,x).$$
We also denote by $\Gamma_{\varepsilon i}$ the functor defined as in the lemma \ref{functor} : 
$$\begin{displaystyle}2\!\!\!\! \varprojlim_{\Pi_2(F_{i})}\alpha \end{displaystyle}\longrightarrow \begin{displaystyle}2\!\!\!\! \varprojlim_{\Pi_2(F_{i+1})}\alpha \end{displaystyle},$$
Using the same process as in the lemma \ref{isomorphism} we define an isomorphism $\Gamma_{\varepsilon}$ : 
$$\xymatrix{
\begin{displaystyle}2\!\!\!\varprojlim_{\Pi_{2}(F_{0})}\alpha \end{displaystyle}\UN[rrr]{\Gamma_{\varepsilon0}}{\Gamma_{\varepsilon1}}{\Gamma_{\varepsilon}}&&&\begin{displaystyle}2\!\!\!\varprojlim_{\Pi_{2}(F_{1})}\alpha \end{displaystyle}.
}$$
~\\

For the rest of this section, for all path $\gamma$ of $B$ and for all homotopy $\varepsilon$ of paths of $B$, we fix a trivialization $h_{\gamma}$ of $I\times_{\gamma}X$ and $h_{\varepsilon}$ of $(I\times I)\times_{\varepsilon}X$ satisfying the condition (\ref{equation_projection}).  \\\\

Let $\gamma, \gamma'$ be two of composable paths in $B$, let $h_{\gamma}$, $h_{\gamma'}$ and $h_{\gamma\land\gamma'}$ the fixed trivializations of  $I\times_{\gamma}X$, $I\times_{\gamma'}X$ and $I\times_{\gamma\land\gamma'}X$ respectively. \\
Let us remark that the map $h_{\gamma}\land h_{\gamma'}$ is also a trivialization of $I\times_{\gamma\land\gamma'}X$. Let us consider the homotopy from $h_{\gamma}\land h_{\gamma'}$ to $h_{\gamma\land\gamma'}$ defined in the lemma \ref{trivialization} and the isomorphism associated by the lemma \ref{isomorphism}. We denote by $H_{\gamma\gamma'}$ the composition of the previous isomorphism with the isomorphism defined by the data of $Id\circ \alpha_{\gamma\gamma'}$ for all $x\in F_{1}$. $H_{\gamma\gamma'}$ is an isomorphism going from $\Gamma_{h_{\gamma}}\circ\Gamma_{h_{\gamma'}}$ to $ \Gamma_{h_{\gamma\gamma'}}$.

We define the $2$-functor $p_{*}$ going from $Rep(\Pi_{2}(X), \Cc\Ac\Tc)$ to \linebreak$Rep(\Pi_{2}(B),\Cc\Ac\Tc)$. 
\begin{Def}
Let $\alpha\in Rep(\Pi_{2}(X), \Cc\Ac\Tc)$, we denote by $p_{*}(\alpha)$ the $2$-functor from $\Pi_{2}(B)$ to $\Cc\Ac\Tc$ that associates  : 
\begin{itemize}
\item to every $b\in B$, the category : 
$$p_{*}(\alpha)(b)=2\!\!\!\! \varprojlim_{\Pi_2(F_{b})}\alpha,$$
\item to every path $\gamma : b_{0}\rightarrow b_{1}$ in $B$, the functor defined in the lemma \ref{functor} :
$$p_{*}(\alpha)(\gamma) =\Gamma : 2\!\!\!\! \varprojlim_{\Pi_2(F_{0})}\alpha \longrightarrow 2\!\!\!\! \varprojlim_{\Pi_2(F_{1})}\alpha,$$
\item to every couple of composable paths $\gamma, \gamma'$ in $B$, the isomorphism of functors $H_{\gamma,\gamma'}$.
$$p_{*}(\alpha)_{\gamma\gamma'} : \xymatrix{\relax\begin{displaystyle} 2\!\!\!\! \varprojlim_{\Pi_2(F_{0})}\alpha \end{displaystyle}\UN[rrr]{}{}{H_{\gamma\gamma' }} & &&\begin{displaystyle}2\!\!\!\! \varprojlim_{\Pi_2(F_{1})}\alpha \end{displaystyle}}
,$$
\item to every homotopy $\varepsilon : \gamma_{0}\rightarrow \gamma_{1}$ in $B$, the composition :
$$ p_*(\varepsilon) : \Gamma_{\gamma_{0}}\buildrel H_{\gamma_{0}}\over\longrightarrow \Gamma_{\varepsilon 0}\buildrel\Gamma_{\varepsilon}\over\longrightarrow \Gamma_{\varepsilon 1}\buildrel H_{\gamma_{1}}\over\longrightarrow \Gamma_{\gamma_{1}},$$
\end{itemize}
The unicity, up to homotopy, of the homotopies between two trivializations assures that these data satisfy the commutation conditions and define a $2$-functor from $\Pi_{2}(X)$ to $\Cc\Ac\Tc$.
And, by the definition of the $2$-limit, this map is $2$-functorial. Hence we denote by $p_{*}$ the $2$-functor from $Rep\big(\Pi_{2}(X),\Cc\Ac\Tc\big)$ to $Rep\big(\Pi_{2}(B),\Cc\Ac\Tc\big)$ defined by these data. 
\end{Def}
\begin{thm}\label{projection}
The following diagram commutes up to equivalence of $2$-functors \nolinebreak:  
$$\xymatrix{
\LLL_{X}\ar[r]^{p_{*}} & \LLL_{B}\ar[d]^{\mu}\\
2Rep\big(\Pi_{2}(X), \Cc\Ac\Tc\big)\ar[r]_{p_{*}}\ar[u]^{\nu} & 2Rep\big(\Pi_{2}(B), \Cc\Ac\Tc\big).
}$$
\end{thm}
\begin{proof}
Let $\alpha$  be a $2$-representation of $\Pi_{2}(B)$ and $\CCC$ the image of $\alpha$ by $\nu$. 
\begin{itemize}
\item Let $b$ a point of $B$, we have : 
$$\big(p_{*}(\CCC)\big)_{b}\simeq\Gamma(F_b,\CCC).$$
Let us recall that by definition $\Gamma(F_b, \CCC)=\Gamma(F_b,i_{F_b}^{-1}\CCC)$, where $i_{F_b}$ is the inclusion of $F_b$ in $X$. Thus, by the proposition \ref{restriction}, we have the natural equivalence :
$$\big(p_{*}(\CCC)\big)_{b}\simeq 2\!\!\!\! \varprojlim_{\Pi_2(F_{b})}\alpha. $$
\item Let $\gamma : b_{0}\rightarrow b_{1}$ be a path in $B$. \\
The proposition \ref{restriction} assures that, for all $x \in X$, there exists natural isomotphism :
$$\xymatrix @!0 @C=1.5cm @R=0.7cm{\Gamma(X,\CCC_{ F_{b}}) \ar[rrr]^{\sim}  \ar[ddd] &&&\begin{displaystyle} 2\!\! \varprojlim_{\Pi_2(F_{b})}\alpha \end{displaystyle} \ar[ddd] \\
 && \ar@{=>}[ld]_{\sim}\\
& & &\\
(\CCC_{F_b})_{x} \ar[rrr]_{\sim} &&& \alpha(x).
}$$
Thus, to show that $\nu(p_*\CCC)(\gamma)\simeq p_*(\alpha)(\gamma)$, it is sufficient to show that the diagram commutes ut to isomorphism : 
$$\xymatrix @!0 @C=1.5cm @R=0.7cm{\Gamma(X,\CCC_{ F_{0}}) \ar[rrr]^{\mu(p_{*}\CCC)(\gamma)}  \ar[ddd] &&&\Gamma(B\times F,\CCC_{ F_{1}}) \ar[ddd] \\
 && \ar@{=>}[ld]_{\sim}\\
& & &\\
\CCC_{x} \ar[rrr]_{\mu(\CCC)\big(\Gamma(\cdot, x)\big)} &&& \CCC_{y}
}$$
where the two vertical functors are the natural restrictions.\\
Using the base-change theorem applied to the first diagram of (\ref{diag}),  we show that, for all $y\in F_{1}$ there exists  isomorphisms of functors visualized by :
$$\xymatrix{
\big(\gamma^{-1}p_{*}\CCC\big)_{0} \eq[d]& \Gamma\big(I, \gamma^{-1}p_{*}\CCC\big)_{1} \ar[r]^{\sim} \ar[l]_{\sim} \eq[d] & \big(\gamma^{-1}p_{*}\CCC\big) \eq[d]\\
\Gamma\big(\{0\}\times F_{0}, \Gamma^{-1}(\CCC)\big) \ar[d]& \Gamma\big(I\times F_{0}, \Gamma^{-1}(\CCC)\big) \ar[d] \ar[l]_{\sim} \ar[r]^{\sim} & \Gamma\big(\{1\}\times F_{0}, \Gamma^{-1}(\CCC)\big) \ar[d]\\
\Gamma^{-1}(\CCC)_{x} & \Gamma(I\times\{y\}, \Gamma^{-1}(\CCC)) \ar[r] \ar[l] & \Gamma^{-1}(\CCC)_{y}
}$$
where $y=\Gamma(1,x)$.\\
This  shows the existence of the isomorphism we looked for.

\end{itemize}
\end{proof}


\section{The $2$-category of stacks on stratified spaces}
Let $(X, \Sigma)$ be a stratified space. It is a natural question to ask if a sheaf is entirely determined by its restrictions on the strata. In other words, if the category of sheaves on $X$ is equivalent to the category whose objects are given by a sheaf on each stratum. The answer is no. To define a sheaf we need some extra data : the gluing data. These are a set of morphisms of sheaves satisfying commutation conditions.

 The following section is a generalization of this problem in the case of stacks. Hence we define a $2$-category whose objects are  the data of a stack on each stratum plus some functors of stacks and morphisms of functors of stacks satisfying some commutation conditions and we show that this category is $2$-equivalent to the $2$-category of stacks on $X$.

Let $\SSS t_{X}$ be the $2$-category of stacks on $X$.
Let us denote $S_{k}$ the union of the strata of dimension $k$ and $i_{k}$ the inclusion of $S_{k}$ on $X$. If $k<l$ we denote by $i_{kl}$ the $2$-functor from $\SSS t_{\Sigma_{l}}$ to $\SSS t_{\Sigma_{k}}$,$i_{kl}=i_{k}^{-1}i_{l*}$. Let us denote by $\eta_{l}$ the $2$-adjunction, $\eta_{l}:Id\longrightarrow i_{l*}i_{l}^{-1}$.\\
For source of simplicity, if $k<l<m$ we denote by $\eta_{l}$ the $2$-functor $Id\bullet \eta_{l}\bullet Id$ going from $i_{km}$ to $i_{kl}i_{lm}$.
 Let us define the $2$-category $\SSS_{\Sigma}$. 
\begin{Def}\label{SSS_{Sigma}}
Let $\SSS_{\Sigma}$ be  the $2$-category defined as follows.
\begin{itemize}
\item[$\bullet$] The objects of $\SSS_{\Sigma}$ are the data  :
\begin{itemize}
\item for every $S_{k}$, a stack $\CCC_{k}$ on $S_{k}$,
\item for every pair $(k,l)$ such that $0\leq k<l\leq n$ (i.e. $S_{k}\subset \overline{S}_{l}$), a functor of stacks : $$F_{kl} : \CCC_{k} \longrightarrow i_{kl}\CCC_{l},$$
\item for every triple $(k,l,m)$ such that $0\leq k<l<m\leq n$ (i.e. $S_{k}\subset\overline{S_{l}}\subset\overline{S_{m}}$), an isomorphism of functors $f_{klm}$ visualized by:
$$
\xymatrix @!0 @C=1.3cm @R=0.4cm {\CCC_{k} \ar[rrr]^{F_{lk}}  \ar[ddddd]_{F_{mk}} &&& i_{kl}\CCC_{l} \ar[ddddd]^{i_{kl}F_{ml}} \\
~\\
& & \ar@{=>}[ld]_\sim^{f_{klm}}\\
&~\\
\\
i_{km}\CCC_{m} \ar[rrr]_{\eta_{l}} &&& i_{kl}i_{lm}\CCC_{m}
}$$ such that the two suitable compositions of the morphisms given by the faces of the following cube : 
$$
\xymatrix {
    \CCC_{k} \ar[rr]^{F_{lk}} \ar[dd]_{F_{pk}} \ar[dr]^{F_{mk}} && i_{kl}\CCC_{l} \ar[dr]^{i_{kl}F_{ml}} \ar[dd]|!{[rd];[ld]}\hole  \\
    & i_{km}\CCC_{m} \ar[rr]^{\eta_{l}} \ar[dd] &\ar[d]& i_{kl}i_{lm}\CCC_{m} \ar[dd]^{i_{kl}i_{lm}F_{pm}} \\
     \ar[rr] |!{[ur];[dr]}\hole^{~~~~~~~~~~~~\eta_{l}} i_{kp}\CCC_{p} \ar[dr]_{\eta_{m}} && i_{kl}i_{lp}\CCC_{p} \ar[rd]^{i_{kl}\eta_{m}} \\
    &  i_{km}i_{mp}\CCC_{p}\ar[rr]_{\eta_{l}} &&  i_{kl}i_{lm}i_{mp}\CCC_{p}\\
  }$$
going from $i_{kl}i_{lm} F_{pm} \circ i_{kl} F_{ml} \circ F_{lk} $ to $i^{-1}_k\eta_l \circ i_{k}^{-1}\eta_{mp} \circ F_{pk}$ are equal.
This means that the following diagram commutes :
$$\xymatrix@!0 @C=3.cm @R=0.9cm {  
                       i_{kl}i_{lm} F_{pm} \circ i_{kl} F_{ml} \circ F_{lk} \ar[rr]^{Id \bullet f_{klm}} \eq[d]^{}&  &  i_{kl}i_{lm} F_{pm} \circ \eta_l \circ F_{mk} \ar[ddd] \\
                       i_{kl}(i_{lm}F_{pm} \circ F_{ml}) \circ F_{lk} \ar[dd]_{i_{kl}f_{lmp}\bullet Id}&  & \\
                       &&\\
                       i_{kl}(\eta_{m} \circ F_{pk}) \circ F_{lk} \eq[d]&&\eta_{l}\circ i_{km}F_{pm}\circ F_{mk}\ar[ddd]^{Id \bullet f_{kmp}} \\
                       i_{kl}\eta_{m} \circ  i_{kl}F_{pl} \circ F_{lk}\ar[dd]_{Id \bullet f_{klp}}  \\
                       \\
                        i_{kl}\eta_{m} \circ \eta_{l} \circ F_{pk} \ar[rr]&  & \eta_l \circ \eta_{m} \circ F_{pk}
                       }$$   
\end{itemize}
\item[$\bullet$] The $1$-morphisms from $\big(\{\CCC_{k}\}, \{F_{kl}\}, \{f_{klm}\}\big)$ to $\big(\{\CCC'_{k}\}, \{F'_{kl}\}, \{f'_{klm}\}\big)$ are given by  : 
\begin{itemize}
\item for every $k\in \{0,\ldots, n\}$, a functor of stacks : $G_{k} : \CCC_{k}\rightarrow \CCC_{k}'$,
\item for every $k,l$ such that $0\leq k<l\leq n$, an isomorphism of functors \nolinebreak:
$$g_{kl} : F'_{lk}\circ G_{k} \buildrel\sim\over\rightarrow i_{kl}G_{l}\circ F_{lk},$$
such that the following diagram commutes :
$$\xymatrix@!0 @C=6cm @R=2cm {  
i_{kl}F'_{ml}\circ F'_{lk}\circ G_{k} \ar[r]^{Id \bullet g_{lk}} \ar[d]_{f_{mlk}\bullet Id} & i_{kl}F_{ml}'\circ i_{kl}G_{l}\circ F_{lk}\ar[d]^{i_{kl}g_{ml} \bullet Id}\\
\eta_{l}\circ F'_{mk}\circ G_{k} \ar[d]_{Id \bullet g_{mk}} & i_{kl}i_{lm}G_{m} \circ i_{kl}F_{ml}\circ F_{lk} \ar[d]^{Id \bullet f_{mlk}}\\
\eta_{l}\circ i_{km}G_{m}\circ F_{mk}\ar[r]&  i_{kl}i_{lm}G_{m} \circ \eta_{l}\circ F_{mk} 
}$$
\end{itemize}
\item[$\bullet$] the $2$-morphisms from the $1$-morphism $\big(\{G_{k}\}, \{g_{kl}\}\big)$ to the $1$-morphism $\big(\{G'_{k}\}, \{g_{kl}\}\big)$ are the data for all $0\leq k\leq n$ of a morphism of functors of stacks $\phi_{k} : G_{k}\rightarrow G'_{k}$, such that the following diagram commutes :
$$\xymatrix@!0 @C=5cm @R=2cm { 
F'_{kl}\circ G_{k} \ar[r]^{g_{kl}} \ar[d]_{Id \bullet \phi_{k}} & i_{kl}G_{l}\circ F_{lk}\ar[d]^{i_{kl}\phi_{l}\bullet Id}\\
F'_{kl}\circ G'_{k} \ar[r]_{g'_{kl}} & i_{kl}G_{l}' \circ F_{lk}
}$$
\end{itemize}
\end{Def}
Hence the objects of this $2$-category are the data of a stack on each stratum plus some gluing data : the functors of stacks and isomorphisms of functors.  \\

To show that $\SSS_{\Sigma}$ is $2$-equivalent to $\SSS t_{X}$ we define two $2$-functors quasi-$2$-inverse to each other : $R_{\Sigma}$ the ``restriction functor''  going from $\SSS t_{X}$ to $\SSS_{\Sigma}$, and $G_{\Sigma}$ the ``gluing functor'' from $\SSS_{\Sigma}$ to $\SSS t_{X}$.
The functor $R_{\Sigma}$ is defined thanks to the restriction and the $2$-adjunction between the $2$-functors $i_{k*}$ and $i_{k}^{-1}$.
\begin{Def}
Let $R_{\Sigma}$ be the $2$-functor going from $\SSS t_{X}$ to $\SSS_{\Sigma}$ which associates to each stack on $X$ the set of its restrictions to each stratum, its adjunction functors and isomorphisms :
$$\begin{array}{cccc}
R_{\Sigma}: & \SSSt _{X}  & \longrightarrow & \SSS_{\Sigma}\\
     & \CCC & \longmapsto &  (\{ \CCC \mid_{S_{k}}\}_{k\leq n }, \{i^{-1}_k\eta_l \}_{k<l\leq n}, \{\lambda_{mlk}\}_{k<l<m\leq n} )\\
     & G : \CCC \rightarrow \CCC' & \longmapsto & ( \{G \mid_{S_{k}}\}_{k \leq n}, \{g_{lk}\}_{k<l\leq n})\\
     & \phi : G \rightarrow G' & \longmapsto & (\{\phi \mid_{S_{k}}\}_{k\leq n})
\end{array}$$
where $\eta_{l}$ is the natural functor 
$$ \eta_{l}: \CCC  \longrightarrow i_{l*}i_{l}^{-1}\CCC,$$
  $\lambda_{mlk}$ are the natural isomorphisms : 
$$\xymatrix @!0 @C=3.2pc @R=2pc {\CCC\mid_{S_{k}} \ar[rrr]  \ar[ddd] &&& (i_{l*}i_{l}^{-1}\CCC)\mid_{S_{k}} \ar[ddd] \\
 && \ar@{=>}[ld]_\sim^{\lambda_{mlk}}\\
& & &\\
(i_{m*}i_{m}^{-1}\CCC)\mid_{S_{k}} \ar[rrr]_{i^{-1}_k\eta_{l}} &&& (i_{l*}i_{lm}i_{m}^{-1}\CCC)\mid_{S_{k}}
}$$
and $g_{lk}$ are the isomorphisms coming from the fact that $\eta_{l}$ is a $2$-transformation \nolinebreak :
$$\xymatrix @!0 @C=3,2pc @R=2pc {\CCC\mid_{S_{k}} \ar[rrr]^{G\mid_{S_{k}}}  \ar[ddd]_{i_k^{-1}\eta_l} &&& \CCC'\mid_{S_{k}} \ar[ddd]^{i_k^{-1}\eta'_l} \\
 && \ar@{=>}[ld]_\sim^{g_{lk}}\\
& & &\\
(i_{l*}i_{l}^{-1}\CCC)\mid_{S_{k}} \ar[rrr]_{(i_{l*}i_{l}^{-1}G)\mid_{S_{k}}} &&& (i_{l*}i_{l}^{-1}\CCC')\mid_{S_{k}}
}$$
As these data come from the adjunction, the commutation conditions are satisfied and the image of $R_{\Sigma}$ belongs to $\SSS_{\Sigma}$.
\end{Def}
Let us define $G_{\Sigma}$. The definition of $G_{\Sigma}$ is inspired by the demonstration of the basic property of gluing stacks on an open covering. Hence the image by $G_{\Sigma}$ of an object $\CCC$ of $\SSS_{\Sigma}$ is a $2$-limit of $i_{k*}\CCC_{k}$, where the $2$-limit encode the gluing data. That is why we define a category $\III$ and, for all object $\CCC$ of $\SSS_{\Sigma}$, a $2$-functor from $\III$ to $\SSS t_{X}$.

\begin{Def}
Let $\III$ be the category defined as follows.
\begin{itemize}
\item Objects of $\III$ are the singletons $\{j\}$ with $0\leq j\leq n$, the couples $(j,k)$ such that $0\leq j<k\leq n$ and the triple $(j,k,l)$ such that $0\leq j<k<l\leq n$. 
\item Morphisms of $\III$ are the data for all objects of $\III$ of :
$$\begin{array}{rccc}
Hom(i,i)&=& \{Id_{i}\}&  \\
Hom((j,k),j)&=&\{s_{jk}^j\}\\
Hom((j,k,l,),(j,k))&=&\{s_{jkl}^{jk}\} &\\
Hom((j,k,l),j)&=&\{s_{jkl}^{j}\} &\\
Hom((j,k,l,),(j,l))&=&\{s_{jkl}^{jl}\}\\
\end{array}$$
\end{itemize}
\end{Def}
 Let $\CCC=\big(\{\CCC_{k}\}, \{F_{kl}\}, \{f_{klm}\}\big)$ be an object of $\SSS_{\Sigma}$, we also denote $\CCC$ the $2$-functor going from $\III$ to $\SSS t_{\Sigma}$ :
$$\begin{array}{cccc}
\CCC : &\III & \longrightarrow & \SSS t_{\Sigma}\\
\end{array}$$
defined as follows :
\begin{itemize}
\item
For every objects $\{k\}$, $\{(k,l)\}$ and $\{(k,l,m)\}$ of $\III$  : 
$$\begin{array}{cclc}
\cat C(j)&=&i_{j*}\CCC_{j}& \\
\cat C (j,k)&=& i_{j*}i_{jk}\CCC_{k} & \\
\cat C (j,k,l)&=& i_{j*}i_{jk}i_{kl}\CCC_{l} & .\\
\end{array}$$
\item
The images of the morphisms are defined as follows : 
\begin{itemize}
\item for every  $s_{jk}^j : (j,k)\rightarrow \{j\}$ such that $j<k$ :
$\CCC(s_{jk}^j) =i_{j*} F_{kj},$
\item for every $s_{jkl}^{jl} : (j,k,l) \rightarrow (j,l)$ :
$\CCC(s_{jkl}^{jl})= \eta_{k}: i_{j*}i_{jl}\CCC_{l}\rightarrow i_{j*}i_{jk}i_{kl}\CCC_{l}.$
\item for every $s_{jkl}^j : (j,k,l) \rightarrow j$ :
$\CCC(s_{jkl}^j) = \eta_{k} \circ i_{j*} F_{lj}.$
\end{itemize}
\item
If $a$ is an object of $\III$ the $2$-morphism $\CCC_{a}: Id_{\CCC(a)} \rightarrow \CCC(Id_{a})$ is the identity. 
\item

If $s$ and $s'$ are two composable morphisms of $\III$ , let us define the $2$-morphism $\CCC_{s,s'}$   :
$$\CCC_{s,s'} : \CCC(s\circ s') \buildrel\sim\over\longrightarrow \CCC(s)\circ \CCC(s')$$
The only two couples of composable morphisms are $(s_{jkl}^{jk}, s_{jk}^j)$ and $(s_{jkl}^{jl},s_{jl}^j)$. We define $\CCC_{s_{jkl}^{jk},s_{jl}^l}$ and $\CCC_{s_{jkl}^{jl},s_{jl}^l}$ as 
$$\CCC_{s_{jkl}^{jk}, s_{jk}^{j}}= f_{lkj},~~\CCC_{s_{jkl}^{jl},s_{jl}^l}=Id$$
\end{itemize}
We define the image of $\CCC$ by $G_{\Sigma}$ by the $2$-limit : $$G_\Sigma(\CCC):=\begin{displaystyle} 2\varprojlim_{\III}\CCC\end{displaystyle}.$$
If $G=\big(\{G_{k}\}, \{g_{kl}\}\big): \CCC\rightarrow \CCC'$ is a $1$-morphism of $\SSS_{\Sigma}$ the commutation conditions satisfied by $g_{kl}$ assure that we can define a functor from the $2$-functor $\CCC$ to $\CCC'$. Taking the $2$-limit we define a functor $G$ from $\begin{displaystyle} 2\varprojlim_{\III}\CCC\end{displaystyle}$ to $\begin{displaystyle} 2\varprojlim_{\III}\CCC'\end{displaystyle}$.\\
In the same way, if $\phi=\{\phi_{k}\} : G\rightarrow G'$ is a $2$-morphism of $\SSS_{\Sigma}$, we can define a morphism between the functors $G$ and $G'$. That is the image of $\phi$ by $G_{\Sigma}$.\\
Now let us consider the $2$-category $2\cat F(\III, \SSS t_{X})$ of $2$-functors from $\III$ to $\SSS t_{X}$. Taking the $2$-limit can be view as a $2$-functor from $\cat F(\III, \SSS t_{X})$ to $\SSS t_{X}$, for a demonstration see for example \cite{St}. Hence we can define the $2$-functor $G_{\Sigma}$.
\begin{Def}
Let $G_{\Sigma}$ be the $2$-functor going from $\SSS_{\Sigma}$ to $\SSS t_{X}$ defined by : 
$$\begin{array}{ccllc}
G_{\Sigma} : & \SSS_{\Sigma} & \longrightarrow & \SSS t_{X}\\\\
& \CCC=\big(\{\CCC_{k}\}, \{F_{kl}\}, \{f_{klm}\}\big) &\longmapsto & \begin{displaystyle} 2\varprojlim_{ \III}\CCC\end{displaystyle}\\
& G: \CCC\rightarrow \CCC' & \longmapsto & \begin{displaystyle} 2\varprojlim_{ \III} \CCC \buildrel G\over\rightarrow 2\varprojlim_{ \III} \CCC'\end{displaystyle}\\
&\phi : G\rightarrow G' &\longmapsto& \xymatrix{2\varprojlim \CCC \UN[r]{}{}{\phi}& 2\varprojlim \CCC' }
\end{array}$$
\end{Def}
\noindent
\textbf{Remarks}
\begin{itemize}
\item[-]
We can define explicitly the stack image of an object of $\SSS_{\Sigma}$. If $U$ is an open of $X$ and $\CCC=\big(\{\CCC_{k}\}, \{F_{kl}\}, \{f_{klm}\}\big)$ an object of $\SSS_{\Sigma}$ is given by :
$$G_{\Sigma}(\CCC)(U)=\big(\{S_{k}\}, \{g_{kl}\}\big)$$
where $S_{k} \in \CCC_{k} (U\cap S_{k})$ and $g_{kl}$ is an isomorphism from $ F_{lk}(S_{k})$ to  $\eta_{kl}(S_{l})$.
\item[-] The commutation conditions satisified by the objects, the $1$-morphims and the $2$-morphisms of $\SSS_\Sigma$ are not necessary to define the functor $G_{\Sigma}$. But without them $G_{\Sigma}$ is not an equivalence.
\end{itemize}
\begin{thm}\label{thm_principal}
The categories $\SSS t_{X}$ and $\SSS_{\Sigma}$ are $2$-equivalent and the functors $R_{\Sigma}$ and $G_{\Sigma}$ are quasi-$2$-inverse. 
\end{thm}
\dem
Let us define two equivalences of $2$-functors  
$$\begin{array}{cccc}
 R_{\Sigma}G_{\Sigma} &\longrightarrow & Id\\
Id & \longrightarrow & G_{\Sigma} R_{\Sigma}.
\end{array}$$
We only define the functor on the objects of the $2$-category, but, as we only use $2$-functor and projection of the $2$-limit to define them, and thanks to the commutations conditions, it is straightforward to show that these applications are $2$-functorial. 

Let us give some notations. For a morphism $s :a\rightarrow b$ of $\III$, we denote by $\pi_{a}$ and $p_{s}$ the projections and the equivalence of functors given by the $2$-limit : 
$$\xymatrix @!0 @C=0.7cm @R=0.7cm {
&&&\begin{displaystyle}2\varprojlim_{c \in \III} \CCC(c) \end{displaystyle} \ar[ddddlll]_{\pi_{a}} \ar[ddddrrr]^{\pi_{b}}\\
~\\
&&&&~\\
&&\ar@{=>}[urr]_{p_{s}}^\sim\\
\CCC(a) \ar[rrrrrr]_{\CCC(s)}&&&&&&\CCC(b)\\
}$$
Let $j$ be an integer smaller than $n$. We also denote $\pi_{a}$ and $p_{s}$, projections and equivalence given by the $2$-limit $\begin{displaystyle} 2\varprojlim_{a \in \III} i_{j}^{-1}\CCC(a)\end{displaystyle}$.
\\
As  $i_{j}$ is an inclusion, the following $2$-natural transform is an equivalence :
$$\varepsilon_{j} : i_{j}^{-1}i_{j*} \longrightarrow Id$$
Let us fix, $(\varepsilon_{j})^{-1}$,  a quasi-inverse of $\varepsilon_{j}$  and $e_{j}$  an isomorphism :
$$e_{j} : \varepsilon_{j} \circ (\varepsilon_{j})^{-1} \longrightarrow Id.$$ 
Let $\Psi_{j}$ denote the functor defined by 
$$\begin{displaystyle}\Psi_{j} : 2\varprojlim_{a \in \III}i_{j}^{-1} \CCC(a) \buildrel\pi_{j}\over\longrightarrow i_{j}^{-1}i_{j*}\CCC_{j} \buildrel\varepsilon_{j}\over\longrightarrow \CCC_{j}\end{displaystyle}$$
This functor an its inverse is essential in the definition of the equivalence between $Id$ and $G_{\Sigma}R_{\Sigma}$.
\begin{lem}
The functor $\Psi_{j} :2\varprojlim i_{j}^{-1} \CCC(a)\rightarrow \CCC_{j} $ is an equivalence. 
\end{lem}
\dem
Let us define $\Phi_{j}$, an inverse of $\Psi_{j}$. \\
As we want to define a functor going to a $2$-limit, it is sufficient to give for all $a$ object of $ \III$, a functor $\Phi_{j}^{a}: \CCC_{j} \rightarrow \CCC(a)$ and for all $s: a\rightarrow b$ morphism of $\III$ an equivalence $h_{s}^j$ of functors :
$$\xymatrix @!0 @C=0.7cm @R=0.7cm {
&&&\CCC_{j} \ar[ddddlll]_{\Phi_{j}^{a}} \ar[ddddrrr]^{\Phi_{j}^b}\\
~\\
&&&&~\\
&&\ar@{=>}[urr]_{h_{j}^{s}}^\sim\\
\CCC(a) \ar[rrrrrr]&&&&&&\CCC(b)\\
}$$
satisfying some commutation conditions. Let us first remark that for $j>k$ and $k<l<m$ we have :
$$\begin{array}{cl}
i_{j}^{-1}\CCC(k)=0 \\
i_{j}^{-1}\CCC(k,l)=0\\
i_{j}^{-1}\CCC(k,l,m)=0.\\
\end{array}$$
Hence we need to define a family of functors for all $j\leq k<l<m$:
\begin{itemize}
\item 
let us recall that $i_{j}^{-1}\CCC(j)= i_{j}^{-1}i_{j*}\CCC_{j}$, we define $\Phi_{j}^j$ by :
$$\Phi_{j}^j=(\varepsilon_{j})^{-1}$$
\item 
for $j< k$ we have $\CCC(k)=i_{jk}\CCC_{k}$, we define $\Phi_{j}^k$ by
$$\Phi_{j}^k=F_{kj}$$
\item 
for $j\leq k<l$ we have $i_{j}^{-1}\CCC(k,l)=i_{jk}i_{kl}\CCC_{l}$, we define $\Phi_{j}^{(k,l)}$ by the composition \nolinebreak:
$$\Phi_{j}^{(k,l)}= \eta_{k}\circ F_{lj}$$
visualized by :
$$ \Phi_{j}^{(k,l)} : \CCC_{j} \buildrel F_{lj}\over\longrightarrow i_{jl}\CCC_{l}   \buildrel\eta_{k}\over\longrightarrow i_{jk}i_{kl}\CCC_{l}$$
\item 
for $j\leq k<l<m$, we have $i_{j}^{-1}\CCC(k,l,m)=i_{jk}i_{kl}i_{lm}\CCC_{m}$, we define $\Phi_{j}^{(k,l,m)}$ by the composition : 
$$\Phi_{j}^{(k,l,m)}=\eta_{l}\circ \eta_{k}\circ F_{mj}.$$
visualized by :
$\CCC_{j} \buildrel F_{mj}\over\longrightarrow i_{jm}\CCC_{m} \buildrel\eta_{k} \over\longrightarrow i_{jk}i_{km}\CCC_{m} \buildrel\eta_{l}\over\longrightarrow i_{jk}i_{kl}i_{lm}\CCC_{m}.$
\end{itemize}
Now let us define the isomorphisms $h_{j}^s$, 

\begin{itemize}
\item for $s=s_{kl}^l : (k,l)\rightarrow l$ we define $h_{j}^s$ by :
$h_{j}^s=Id$
\item for $s=s_{kl}^k :(k,l)\rightarrow k$ with $j<k$ we define $h_{j}^s$ by : 
$h_{j}^s=\theta_{jkl}$
\item let us consider the morphism $s=s_{j}^{jk}: (j,k)\rightarrow j$, the morphism $h_{j}^s$ is going from $i_{j}^{-1}i_{j*}F_{kj} \circ \varepsilon_{j}^{-1}$ to $ \eta_{j}\circ F_{kj}$ :
$$h_{j}^s :  i_{j}^{-1}i_{j*}F_{kj} \circ \varepsilon_{j}^{-1} \buildrel\sim\over\longrightarrow \eta_{j}\circ F_{kj}$$
As $\varepsilon_{j}$ comes from a natural $2$-transform, there exists an isomorphism $\theta_{j} : F_{kj} \circ \varepsilon_{j} \buildrel\sim\over\rightarrow \varepsilon_{j}\circ (i_{j}^{-1}i_{j*}F_{kj})$ :
$$\xymatrix @!0 @C=1cm @R=0,5cm {\CCC_{j}  \ar[ddddd]_{F_{kj}} &&& i_{j}^{-1}i_{j*}\CCC_{j} \ar[lll]_{\varepsilon_{j}} \ar[ddddd]^{i_{j}^{-1}i_{j*}F_{kj}} \\
~\\
& & \ar@{=>}[ld]_{\sim}^{\theta_{j}}\\
&~\\
\\
i_{jk}\CCC_{k} &&& i_{j}^{-1}i_{j*}i_{jk}\CCC_{k}\ar[lll]^{\varepsilon_{j}}.}$$
Hence we have the isomorphism :
$$
\xymatrix@!0 @C=5cm {
\eta_{j}\circ F_{kj}\circ \varepsilon_{j}\circ (\varepsilon_{j})^{-1} \ar[r]^{ ~~~~~~~~~Id~\bullet \theta_{j}\bullet Id}&\eta_{j}}\circ \varepsilon_{j}\circ i_{j*}i_{j}^{-1}F_{kj}\circ (\varepsilon_{j})^{-1}.$$
Now, $\eta_{j}$ is a left quasi-inverse of $\varepsilon_{j}$, let us denote  $n_{j}$ the isomorphism \nolinebreak:
$$n_{j} : \eta_{k}\circ \varepsilon_{j} \longrightarrow Id $$
In the same way let us recall that the isomorphism $e_{j}$ goes from $\varepsilon_{j}\circ (\varepsilon_{j})^{-1}$ to $Id$.
We define $h_{j}^s$ by :
$$h_{j}^s = n_{j} \circ (I_{i_{j}^{-1}\eta_{j}} \bullet \theta_{j} \bullet I_{\varepsilon_{j}^{-1}}) \circ e_{j}$$
\item for $s=s_{klm}^{km}$, we define $h_{j}^s$ by the identity,
\item let us consider the morphism $s_{klm}^{kl}: (k,l,m) \rightarrow (k,l)$, the morphism $h_{j}^s$ is going from the composition $(i_{jk}i_{kl}F_{ml} )\circ \eta_{k} \circ F_{lj}$ to $ \eta_{l} \circ \nolinebreak \eta_{k} \circ F_{mj}$ : 
$$h_{j}^s :(i_{jk}i_{kl}F_{ml} )\circ \eta_{k} \circ F_{lj} \buildrel\sim\over\longrightarrow\eta_{l} \circ \eta_{k} \circ F_{mj}.$$
But, by definition of $\theta_{jlm}$ and as $\eta_{kl}$ comes from a $2$-adjunction, we have the two following isomorphisms :
$$\xymatrix @!0 @C=3,3cm @R=1cm {
& i_{jl}\CCC_{l} \ar[r]^{\eta_{k}} \ar[dddr]_{i_{jl}F_{ml}} &  i_{jk}i_{kl}\CCC_{l} \ar[dddr]^{i_{jk}i_{kl}F_{ml}} \\
&\ar@{=>}[d]^\sim_{\theta_{jlm}} & \ar@{=>}[d]^\sim_{i_{j}^{-1}adj}\\
&~& ~\\
\CCC_{j} \ar[uuur]^{F_{lj}} \ar[r]_{F_{mj}} & i_{jm}\CCC_{m} \ar[r]_{ \eta_{l}}& i_{jl}i_{lm}\CCC_{m} \ar[r]_{\eta_{k}} &  i_{jk}i_{kl}i_{lm}\CCC_{m}.
}$$
Then, $h_{j}^s$ is defined by a the correct composition of these two isomorphisms. 
\end{itemize}
Thanks to the commutation conditions and as the morphisms of adjunction satisfy good conditions of commutation, these functors and isomorphisms of functors satisfy the condition to define a functor coming from $\CCC_{j}$ to $2\varprojlim i_{j}^{-1} \CCC(a)$,  let us denote this functor $\Phi_{j}$ : 
$$\Phi_{j} : \CCC_{j} \longrightarrow 2\varprojlim i_{j}^{-1} \CCC(a)$$
and for every objects $a$ of $\III$, let us denote $\varphi_{j}^{a}$ the isomorphism :
$$\varphi_{j}^a : \pi_{a} \circ \Phi_{j} \buildrel\sim\over\longrightarrow \Phi_{j}^{a}.$$
Moreover, the definition of this functor is $2$-functorial. This is comes from the fact that  the conditions to be an $1$-morphism in the $2$-category $\SSS_{\Sigma}$ and the fact that the $2$-limit can be viewed as a $2$-functor coming from the $2$-functor going from $\III$ to $\SSS t_{X}$ to $\SSS t_{X}$.

It remains to show that $\Phi_{j}$ is a quasi-inverse of $\Psi_{j}$. The easiest part is to show that $\Psi_{j}\circ \Phi_{j}$ is isomorphic to the identity. We have by definition  : 
$$\xymatrix @!0 @C=1.5cm{
\Psi_{j} \circ \Phi_{j} = \varepsilon_{j} \circ \pi_{j} \circ \Phi_{j} \ar[rrr]^{~~~~\sim}_{~~~~I_{\varepsilon_{j}}\bullet \varphi_{j}^j}&&& \varepsilon_{j} \circ (\varepsilon_{j})^{-1}  \ar[rr]_{~~~~e_{j}}^{~~~~\sim} && Id
.}$$
Hence the composition  $ e_{j}\circ (I_{\varepsilon_{j}}\bullet \varphi_{j}^j)$ is an isomorphism between $\Psi_{j}\circ \Phi_{j}$ and $Id$.

Let us show that $\Phi_{j}\circ \Psi_{j}$ is isomorphic to the identity. To show that, it is sufficient to show that for every object $a$ of $\III$ there exists compatible morphisms $l_{j}^{a}$ :
$$l_{j}^{a} : \pi_{a}\circ \Phi_{j} \circ \Psi_{j} \buildrel\sim\over\longrightarrow \pi_{a}.$$

Let $\{k\}$ be an object of $\III$, by definition of $\Phi_{j}$ we have the following isomorphism \nolinebreak :
$$
\varphi_{j}^k : \pi_{k}\circ \Phi_{j} \buildrel\sim\over\longrightarrow F_{mj},
$$
then let us consider the following isomorphism :
\begin{equation}\label{morphism1}
\varphi_{j}^k\bullet I_{\Psi_{j}} : \pi_{k}\circ \Phi_{j}\circ \Psi_{j}\buildrel\sim\over\longrightarrow F_{mj}\circ \Psi_{j}
\end{equation}
where $I_{\Psi_{j}}$ is the morphism identity of $\Psi_{j}$.\\
Then we have the following isomorphisms :
$$\xymatrix @!0 @C=1,5cm @R=1.2cm { 
&&&\ar@{=}[d]\\
 2\varprojlim i_{j}^{-1} \CCC(a) \ar[rrr]^{\pi_{j}} \ar[ddd]_{\pi_{m}} \ar@/^1.5cm/[rrrrrr]^{\Psi_{j}} &&& i_{j}^{-1}i_{j*}\CCC_{j} \ar[ddd]_{i_{j}^{-1}i_{j*}F_{mj}}\ar[rrr]^{\varepsilon_{j}} &&& \CCC_{j} \ar[ddd]_{F_{kj}}\\
&&\ar@{=>}[dl]_{\sim}^{p_{jk}}&&&\ar@{=>}[dl]_{\sim}\\
&&&&&&&\\
i_{jk}\CCC_{k} \ar[rrr]_{i_{j}^{-1}\eta_{j}} \ar@/_1.5cm/[rrrrrr]_{Id}&&& i_{j}^{-1}i_{j*}i_{jk}\CCC_{k}  \ar[rrr]_{\varepsilon_{j}} &&& i_{jk}\CCC_{k}\\
&&& \ar@{=>}[u]_{\sim}
}.$$
Let us recall that $p_{jk}$ is the isomorphism coming from the $2$-limit, the equality is given by the definition of $\Psi_{j}$ and the two others isomorphisms come from the $2$-adjunction. Hence, by composing the morphism above, we can define an isomorphism between $F_{mj}\circ \Psi_{j}$ and $\pi_{m}$ : 
\begin{equation}\label{morphism2}
F_{mj}\circ \Psi_{j} \buildrel\sim\over\longrightarrow \pi_{m}.
\end{equation}
We define $l_{j}^k$ by the vertical composition of the isomorphisms (\ref{morphism1}) and (\ref{morphism2}).

Let $(k,l)$ and $(k,l,m)$ be two objects of $\III$, the $2$-limit give us these isomorphisms : 
$$\xymatrix @!0 @C=1cm @R=1.5cm {
&&&&&2\varprojlim \CCC(a) \ar[dd]_{\pi_{ml}}\ar[ddlllll]_{\pi_{m}} \ar[ddrrrrr]^{\pi_{klm}}\\
&&&\ar@{=>}[r]^\sim_{p_{lm}}&&&\ar@{=>}[r]^\sim_{p_{klm}}&&\\
i_{jm}\CCC_{m}  \ar[rrrr]_{~~~~~~\eta_{l}} &&&&&~~~~i_{jl}i_{lm}\CCC_{m}\ar[rrrr]_{~~~~~~~~~\eta_{k}}&&&&&i_{jk}i_{kl}i_{lm}\CCC_{m}.
}$$
Hence by composing horizontally $p_{lm}$ with the identity of $\Phi_{j}$ and $\Psi_{j}$ we obtain \nolinebreak:
\begin{equation}
\pi_{lm} \circ \Phi_{j} \circ \Psi_{j} \buildrel\sim\over\longrightarrow i_{j}^{-1}\eta_{lm} \circ \pi_{m} \circ \Phi_{j} \circ \Psi_{j}.
\end{equation}
In the same way, by composing horizontally the identity morphism of $\eta_{l}$ and the isomorphism $l_{j}^m$ we obtain the isomorphism :
\begin{equation}
 i_{j}^{-1}\eta_{lm} \circ \pi_{m} \circ \Phi_{j} \circ \Psi_{j} \buildrel\sim\over\longrightarrow  i_{j}^{-1}\eta_{lm} \circ \pi_{m}
\end{equation}
We define $l_{j}^{kl}$ as the vertical composition of the isomorphism (3), (4) and the inverse of $p_{lm}$. The isomorphism $l_{j}^{klm}$ is defined in the same way.

The $2$-functorial feature of the isomorphisms of adjunction and the compatibility of the isomorphisms of projection assure that the commutation conditions are satisfied. Hence, they define an isomorphisms of functors :
$$l_{j} : \Phi_{j} \circ \Psi_{j} \buildrel\sim\over\longrightarrow Id.$$
\cqfd
Let us come back to the theorem \ref{thm_principal}, and let us define an equivalence of $2$-functors :
$$Id \buildrel\sim\over\longrightarrow R_{\Sigma}\circ G_{\Sigma}.$$
We are going to define this equivalence only on objects. The natural feature of the equivalences considered, and the conditions to be an object, a $1$-morphism or a $2$-morphism of $\SSS_{\Sigma}$, assure that the map that we are going to define can be extended in a natural $2$-transform between the $2$-functors $Id$ and $R_{\Sigma}\circ G_{\Sigma}$.

Let $\CCC=\big(\{\CCC_{k}\}, \{F_{kl}\}, \{f_{klm}\}\big)$ be an object of $\SSS_{\Sigma}$. We need to define, for all $j\leq n$, a natural equivalence 
$$\alpha_{j} : \CCC_{j} \buildrel\sim\over\longrightarrow i_{j}^{-1}2\varprojlim \CCC(a)$$
such that there exists for all $j <k\leq n$ an isomorphism :
$$\xymatrix @!0 @C=1.5cm @R=0,6cm {
\CCC_{j} \ar[ddddd] \ar[rrr]^{F_{kj}} &&& i_{jk}\CCC_{k} \ar[ddddd] \\
~\\
& & \ar@{=>}[ld]_{\sim}^{}\\
&~\\
\\
i_{j}^{-1}2\varprojlim \CCC(a)\ar[rrr]_{\eta_{k}}  &&& i_{jk}i_{k}^{-1}2\varprojlim\CCC(a).
}$$

The finite $2$-limits commute up to isomorphism with the inductive $2$-colimits. It is a particular case of a theorem shown in \cite{Dubuc}, for a immediate proof see \cite{Interchange}. Hence  the $2$-functors $i_{k}^{-1}$ and $i_{k*}$  commute up to equivalence with the finite \linebreak$2$-limits. Hence we have the following natural equivalences :
$$2\varprojlim i_{j}^{-1}\CCC(a) \buildrel\sim\over\longrightarrow   i_{j}^{-1}2\varprojlim \CCC(a)$$
$$2\varprojlim i_{jk}i_{k}^{-1}\CCC(a)   \buildrel\sim\over\longrightarrow i_{jk}i_{k}^{-1}2\varprojlim\CCC(a)$$
and the isomorphism :
$$\xymatrix @!0 @C=1.5cm @R=0,6cm {2\varprojlim i_{j}^{-1} \CCC(a)\ar[rrr]^{2\varprojlim \eta_{k}}  \ar[ddddd]_{}  &&&2\varprojlim i_{jk}i_{k}^{-1}\CCC(a)  \ar[ddddd]_{} \\
~\\
& & \ar@{=>}[ld]_{\sim}^{}\\
&~\\
\\
 i_{j}^{-1}2\varprojlim \CCC(a)\ar[rrr]_{\eta_{k}} &&& i_{jk}i_{k}^{-1}2\varprojlim\CCC(a) .}$$
We define $\alpha_{j}$ by  the composition :
$$\alpha_{j} :  \CCC_{j} \buildrel \Phi_{j}\over\longrightarrow 2\varprojlim  i_{j}^{-1}\CCC(a) \buildrel\sim\over\longrightarrow i_{j}^{-1}2\varprojlim \CCC(a)$$
It remains to define a natural equivalence :
$$\xymatrix @!0 @C=1.5cm @R=0,6cm {
\CCC_{j} \ar[ddddd]_{\Phi_{j}} \ar[rrr]^{F_{kj}} &&& i_{jk}\CCC_{k} \ar[ddddd]^{i_{jk}\Phi_{k}} \\
~\\
& & \ar@{=>}[ld]_{\sim}^{}\\
&~\\
\\
2\varprojlim i_{j}^{-1} \CCC(a)\ar[rrr]_{2\varprojlim i_{j}^{-1}\eta_{k}}  &&& 2\varprojlim i_{jk}i_{k}^{-1}\CCC(a) .
}$$
As $\Psi_{k}$ is a quasi-inverse of $\Phi_{k}$ it is sufficient to define an isomorphism between the functors :
$$i_{jk}\Psi_{j} \circ i_{j}^{-1}\eta_{k}\circ \Phi_{j} \buildrel\sim\over\longrightarrow F_{kj}.$$
Now, we have the following isomorphisms :
$$\xymatrix @!0 @C=1.10cm @R=0.5cm { 
 &&&2\varprojlim i_{j}^{-1}\CCC(a)\ar[rrr] \ar[ddddd]_{\pi_{k}} &&&2 \varprojlim i_{jk}i_{k}^{-1}\CCC(a)\ar[ddddd]^{\pi_{k,}}\ar[dddddddddrrr]^{i_{jk}\Psi_{k}}\\
& \\
&&&&& \ar@{=>}[ld]_{\sim}^{}\\
&&&&&&&&&\\
&&\ar@{=>}[dd]_{\sim}&&&&& \ar@{=>}[dd]^\sim\\
&&& i_{j}^{-1}\CCC(k) \ar[rrr]_{i_{j}^{-1}\eta_{k}} &&& i_{jk}i_{k}^{-1}\CCC(k)\ar[ddddrrr]_{i_{jk}\varepsilon_{k}}\\
&&&&&\ar@{=>}[dd]_{\sim}&&&&\\
&&&&\\
&&&&&&&&&\\
\CCC_{j}\ar[uuuuuuuuurrr]^{\Phi_{j}} \ar[rrrrrrrrr]_{F_{kj}}\ar[rrruuuu]_{F_{kj}} &&&&&&&&& i_{jk}\CCC_{k}
.}$$
The two isomorphisms of triangle are given by the definition of $\Phi_{j}$ and $\Psi_{k}$. The isomorphism of the top is given by the $2$-limit and the last one is a horizontal composition of the identity and the isomorphism given by the $2$-adjunction :
$$i_{j}^{-1}\eta_{k} \circ i_{jk}\varepsilon_{k} \simeq Id.$$
The suitable composition of these isomorphisms gives the isomorphism looked for.\\
Hence we have defined an isomorphism from $\CCC$ to $R_{\Sigma}\circ G_{\Sigma}(\CCC)$. This isomorphism is $2$-functorial and this shows that $Id$ is equivalent to $R_{\Sigma}\circ G_{\Sigma}.$\\

Let us define an isomorphism $\beta$ from $Id$ to $G_{\Sigma}\circ R_{\Sigma}$. Let  $\GGG$ be a stack on $X$. Let us recall that from the definition of $R_{\Sigma}$ we have :
$$R_{\Sigma}(\GGG)=(\{i_{k}^{-1} \GGG\}, \{i_{k}^{-1}\eta_{l}\}, \{\lambda_{klm}\})$$
where $\eta_{l}$ is the natural functor :
$$\eta_{l} : \GGG \longrightarrow i_{l*}i_{l}^{-1}\GGG$$
and $\lambda_{klm}$ is the natural isomorphism :
$$
\xymatrix @!0 @C=3.2pc @R=2pc {i_{k}^{-1}\GGG \ar[rrr]^{\eta_{l}}  \ar[ddd]_{\eta_{m}} &&& i_{kl}i_{l}^{-1}\GGG \ar[ddd]^{i_{kl}i_{l}^{-1}\eta_{m}} \\
 && \ar@{=>}[ld]_\sim^{\lambda_{klm}}\\
& & &\\
 i_{km}i_{m}^{-1}\GGG \ar[rrr]_{i^{-1}_k\eta_{lm}} &&&  i_{kl}i_{lm}i_{m}^{-1}\GGG.
}$$
Let us consider the family of functors $\{\Xi_{a}\}_{a\in \III}$ :
\begin{itemize}
\item for every $k\leq n$, the functor $\Xi_{k}$ is defined by :
$$\Xi_{k}=\eta_{k}$$
\item for every pair $(k,l)$ with $k<l\leq n$ let us define $\Xi_{(k,l)}$ by :
$$\Xi_{(k,l)} : \GGG \buildrel\eta_{k}\over\longrightarrow i_{k*}i_{k}^{-1}\GGG \buildrel i_{k*}i_{k}^{-1}\eta_{l}\over\longrightarrow i_{k*}i_{kl}i_{l}^{-1}\GGG$$
\item for every triple $(k,l,m)$ such that $k<l<m\leq n$, $\Xi_{(k,l,m)}$ is defined by : 
$$\xymatrix @!0 @C=1.5cm @R=2cm {
\GGG \ar[r]^{\eta_{k}} & i_{k*}i_{k}^{-1}\GGG \ar[rr]^{i_{k*}i_{k}^{-1}\eta_{l}} && i_{k*}i_{kl}i_{l}^{-1}\GGG \ar[rrr]^{ i_{k*}i_{kl}i_{l}\eta_{m}} &&& i_{k*}i_{kl}i_{l}i_{m*}i_{m}^{-1}\GGG
}$$
\end{itemize}
and the family of isomorphisms of functors $\{\xi_{s}\}_{s \in \Mc or(\III)}$ : 
\begin{itemize}
\item for the morphisms $(k,l)\rightarrow k$ and $(k,l,m)\rightarrow (k,l)$ the isomorphism $\xi_{s}$ is the identity.
\item for the morphism $(k,l) \rightarrow l$ the isomorphism $\xi_{s}$ is the morphism $\lambda_{kl}$ : 
$$\xymatrix @!0 @C=3.2pc @R=2pc{\GGG \ar[rrr]  \ar[ddd] &&& i_{k}^{-1}\GGG \ar[ddd] \\
 && \ar@{=>}[ld]_{\sim}\\
& & &\\
 i_{l*}i_{l}^{-1}\GGG \ar[rrr]_{\eta_{kl}} &&&  i_{k*}i_{kl}i_{l}^{-1}\GGG
}$$
\item for the morphism $(k,l,m)\rightarrow (k,m)$ the functor $\xi_{s}$ is defined by the horizontal composition of the identity, the morphism $\eta_{k}$ and the isomorphism $\lambda_{klm}$. 
\end{itemize}
These families satisfy the compatibility relations. Hence these data define a morphism $\Xi$ :
$$\Xi : \GGG \longrightarrow R_{\Sigma}\circ G_{\Sigma}(\GGG)$$
and for all object $a$ of $\III$ a unique isomorphism $\xi_{a}$ : 
$$\xi_{a} : \Xi\circ \pi_{a}\buildrel\sim\over\longrightarrow \Xi_{a}$$
But the restriction to each stratum is an equivalence. This show that $\Xi$ is an equivalence of functor. 
\cqfd
\section{Constructible stacks}

In this section we consider the $2$-category, $\SSS t^c_{\Sigma}$ of constructible stacks relatively to a fixed stratification $\Sigma$. This notion was introduced by D. Treumann in \cite{Tr1}. It is a natural generalization of the notion of constructible sheaf. 
\begin{Def}
A stack on $X$ is called constructible relatively to $\Sigma$ if its restrictions to each stratum is locally constant. 
We denote by $\SSS t_{\Sigma}^c$ the thick sub-$2$-category of $\SSS t_{X}$ whose objects are locally constant stacks.
\end{Def}

An important example of constructible stack is the stack of perverse sheaves.\\

The aim of this section is to describe the $2$-category $\SSS t_{\Sigma}^c$ in the language of $2$-representations. Here we need to consider more particular stratified spaces : the Thom-Mather spaces relatively $2$-connected. A Thom-Mather space $(X,\Sigma)$ is given, for all strata $\Sigma_{k}$, with a tubular neighborhood $T_{k}$, a locally trivial fiber bundle $p_{k}$ :
$$p_{k}: T_{k}\longrightarrow \Sigma_{k}.$$ 
and a continuous map $\rho_{k}$ from $T_{k}$ to $\RR^+$, named distance map, such that $$\Sigma_{k}=\rho_{k}^{-1}(0).$$
Let us recall that if $\Sigma_{k}\subset \overline{\Sigma}_{l} $, then we have the equality :
$$p_{k}\circ p_{l}\mid_{T_{k}\cap T_{l}}=p_{k}\mid_{T_{k}\cap T_{l}}.$$
Let us denote by $i_{k}$ the natural inclusion of the stratum $\Sigma_{k}$ in $X$. If $\Sigma_{k}\subset \overline{\Sigma}_{l}\subset \overline{\Sigma}_{m}\subset\overline{\Sigma}_{p}$ and if $i_{k,l}$, $i_{k,lm}$ and $i_{k,lmp}$ are the inclusions :
$$i_{l,k} : \Sigma_{l}\cap T_{k}\hookrightarrow \Sigma_{l},~~~i_{m,lk} : \Sigma_{m}\cap T_{l}\cap T_{k}\hookrightarrow \Sigma_{m},~~~i_{p,mlk} : \Sigma_{p}\cap T_{m}\cap T_{l}\cap T_{k}\hookrightarrow \Sigma_{p},$$
we denote by $p_{k,l}$, $p_{k,lm}$, $p_{k,lmp}$, the $2$-functors : 
$$\begin{array}{lclcc}
p_{k,l}=p_{k*}\circ i_{l,k}^{-1} &:&\LLL_{\Sigma_{l}}&\longrightarrow &\LLL_{\Sigma_{k}}\\
 p_{k,lm}=p_{k*}\circ i_{m,lk}^{-1}&:&\LLL_{\Sigma_{m}}&\longrightarrow &\LLL_{\Sigma_{k}}\\
  p_{k,lmp}=p_{k*}\circ i_{p,mlk}^{-1}&:&\LLL_{\Sigma_{p}}&\longrightarrow &\LLL_{\Sigma_{k}}
\end{array}$$
If we denote by $j_{m,lk}$ the inclusion of $\Sigma_{m}\cap T_{l}\cap T_{k}$ in $\Sigma_{m}\cap T_{k}$ and $\eta_{l}$ the $2$-adjunction from $Id$ to $j_{m,lk*}j_{m,lk}^{-1}$, we have that $Id_{p_{k*}}\bullet \eta_{l}\bullet Id_{i_{l,k}^{-1}}$ is going from $p_{k,m}$ to $p_{k,lm}$, we also denote it $\eta_{l}$. 
$$\eta_{l} : p_{k,m}\longrightarrow p_{k,lm}$$
With another harmless abuse of notation we denote by $\eta_{l}$ the natural transformation from $p_{k,mp}$ to $p_{k,lmp}$ defined with the adjunction.

Let $\CCC$ be a constructible stack relatively to $\Sigma$, by definition, $R_{\Sigma}(\CCC)$ is the data of locally constant stacks on each stratum plus the gluing conditions. As we have seen in the first section, a locally constant stack is nothing but a $2$-representation. Hence it remains to express the gluing conditions in terms of $2$-representations. The first thing to do is to verify that the image of a locally constant stack through the functor $i_{k}^{-1}i_{l*}$ is still a locally constant stack. Then we have to express this functor in term of $2$-representations. In order to do this, we show (corollary \ref{Cor}) the equivalence of the $2$-functors $i_{k,l}$ and $p_{k,l}$ restricted to the $2$-category $\SSS t_{\Sigma}^c$. Finally,  the translation of the $2$-functors $i_{k}^{-1}$ and $p_{*}$ in the language of $2$-representations that we have done in the first section allows us to define a combinatoric $2$-category $2$-equivalent to the $2$-category of constructible stacks. 

To do this, the base-change theorem showed by D. Treumann in \cite{Tr1} and the following lemma are very convenient.
\begin{lem}\label{open}
Let $Y$ be a topological space, $V$ be an open of $X$ and $F$ be a subset of $X$. we denote by  $i_{V}$, $i_{F}$, $j_{V}$ and $j_{F}$ the following inclusions :
$$\xymatrix{
V \hook[r]^{i_{V}} &Y\\
V\cap F \hook[u]^{j_{F}} \hook[r]_{j_{V}} & F \hook[u]_{i_{F}}.
}$$
Then the base-change map : 
$$i_{V}^{-1}i_{F*}\buildrel\sim\over\longrightarrow j_{F*}j_{V}^{-1}$$
is an equivalence.
\end{lem}
\dem
It is sufficient to see that the base-change map is an equivalence on the stalks.
\cqfd
\begin{lem}
Let $\Sigma_{k}$, $\Sigma_{l}$, $\Sigma_{m}$ and $\Sigma_{p}$ be four strata such that $\Sigma_{k}\subset \overline{\Sigma}_{l}\subset \overline{\Sigma}_{m}\subset \overline{\Sigma}_{p}$. There exists  natural equivalences of $2$-functors :
$$\begin{array}{cclc}
p_{k,l}\circ p_{l,m}&\buildrel\sim\over\longrightarrow &p_{k,lm}\\
p_{k,l}\circ p_{l,mp}&\buildrel\sim\over\longrightarrow & p_{k,lmp}\\
p_{k,lm}\circ p_{m,p}&\buildrel\sim\over\longrightarrow & p_{k,lmp}.
\end{array}$$
\end{lem}
\dem 
Considering the following commutative diagram :
$$\xymatrix{
\Sigma_{p}\cap T_{m}\cap T_{l}\cap T_{k}\hook[r] \ar[d]_{p_{m}} & \Sigma_{p}\cap T_{m}\cap T_{l} \hook[r] \ar[d]_{p_{m}} &\Sigma_{p}\cap T_{m}\hook[r] \ar[d]_{p_{m}}& \Sigma_{p} \\
\Sigma_{m}\cap T_{l}\cap T_{k} \hook[r] \ar[d]_{p_{l}} & \Sigma_{m} \cap T_{l} \hook[r] \ar[d]_{p_{l}} & \Sigma_{m}\\
\Sigma_{l} \cap T_{k}\ar[d] \hook[r] \ar[d]^{p_{k}}& \Sigma_{l}\\
\Sigma_{k}
}$$
and as $p_{k}\circ p_{l}|_{T_{k}\cap T_{l}}=p_{k}|_{T_{k}\cap T_{l}}$, the lemma is a direct application of the lemma \ref{open}.
\cqfd
\begin{prop}\label{2equivalence}
Let $\Sigma_{k}$ be a stratum and $T_{k}$ be the tubular neighborhood, the $2$-functor $p_{k*}$ restricted to the $2$-category, $\SSSt_{\Sigma}^c$, of constructible stacks on $T_{k}$ goes to the $2$-category $\LLL_{\Sigma_{k}}$. Moreover $p_{k*}$ is equivalent to the $2$-functor $i_{k}^{-1}$.
\end{prop}
\dem
Let us consider the $2$-natural functor $\eta_{k}$ given by the $2$-adjunction : 
$$\eta_{k} : p_{k*}\longrightarrow p_{k*}i_{k*}i_{k}^{-1}.$$
By definition of a Thom-Mather space $p_{k}\circ i_{k}=Id$, thus $p_{k*}i_{k*}i_{k}^{-1}$ is naturally isomorphic to $i_{k}^{-1}$.
Hence, there exists a natural functor from $p_{k*}$  to $i_{k}^{-1}$. Let $\CCC$ be a constructible stack on $T_{k}$, let us show that the functor applied to $\CCC$ is an equivalence on its stalks. Let $x\in \Sigma_{k}$. We denote by $F_{k}$ the set $p_{k}^{-1}(x)$ and $F_{k\varepsilon}=F_{k}\cap \rho_{k}^{-1}([0, \varepsilon[)$.  The family $U\times F_{k\varepsilon}$ is a base of neighborhoods at $x$, hence we have the equivalences \nolinebreak: 
$$\begin{array}{cclc}
(p_{k*}\CCC)_{x} &\simeq &\begin{displaystyle}  2\!\!\!\!\varinjlim_{\substack{x\in U\subset \Sigma_{k} }}\Gamma(U\times F_{k}, \CCC) \end{displaystyle}\\
(i_{k}^{-1}\CCC)_{x} &\simeq &\begin{displaystyle}  2\!\!\!\!\varinjlim_{\substack{x\in U\subset \Sigma_{k}\\ \varepsilon>0 }}\Gamma(U\times F_{k\varepsilon}, \CCC) \end{displaystyle} \\
&\simeq&\begin{displaystyle}  2\!\!\!\!\varinjlim_{\substack{x\in U\subset \Sigma_{k}}}2\varinjlim_{\varepsilon>0}\Gamma(U\times F_{k\varepsilon}, \CCC) \end{displaystyle} 

\end{array}$$
Now, the inclusion $U\times F_{k\varepsilon}\hookrightarrow U\times F_{k}$ is a stratified homotopy equivalence. As $\CCC$ is a constructible stack, this is shown, by proposition 3.13 of \cite{Tr1}, that the $2$-limit $\begin{displaystyle}  2\varinjlim_{ \varepsilon>0 }\Gamma(U\times F_{k\varepsilon}, \CCC) \end{displaystyle} $ is constant an equal to $\Gamma(U\times F_{x}, \CCC)$.
\cqfd
\begin{cor}\label{Cor}
Let $\Sigma_{k}$, $\Sigma_{l}$, $\Sigma_{m}$ and $\Sigma_{p}$ be four strata such that $\Sigma_{k}\subset \overline{\Sigma}_{l}\subset \overline{\Sigma}_{m}\subset \overline{\Sigma}_{p}$, the $2$-functors $p_{k,l}$, $p_{k,lm}$, $p_{k,lmp}$ are respectively equivalent to the $2$-functors $i_{kl}$, $i_{kl}i_{lm}$ and $i_{kl}i_{lm}i_{mp}$. Moreover, there exists isomorphisms of functors :
$$
\xymatrix @!0 @C=1.3cm @R=0.4cm {i_{kp} \ar[rrr]  \ar[ddddd]_{\eta_{l}} &&& p_{k,p} \ar[ddddd]^{\eta_{l}} && i_{kl}i_{lp}\ar[rrr]  \ar[ddddd]_{\eta_{m}} &&& p_{k,lp} \ar[ddddd]^{\eta_{m}}\\
&&&&&&&&~\\
& & \ar@{=>}[ld]_\sim^{}&&&& & \ar@{=>}[ld]_\sim^{}\\
&&&&&&&&~\\
\\
i_{kl}i_{lp} \ar[rrr] &&& p_{k,lp}&&i_{kl}i_{lm}i_{mp} \ar[rrr] &&& p_{k,lmp}
}$$
such that the two suitable compositions of the isomorphisms given by the faces of this cube are equal : 
$$
\xymatrix {
    i_{kp} \ar[rr] \ar[dd] \ar[dr] && p_{k,l} \ar[dr] \ar[dd]|!{[rd];[ld]}\hole  \\
    & i_{km}i_{mp} \ar[rr] \ar[dd] &\ar[d]& p_{k,mp} \ar[dd]^{} \\
     i_{kl}i_{lp} \ar[rr] |!{[ur];[dr]}\hole^{}  \ar[dr]_{} && p_{k,lp} \ar[rd] \\
    &  i_{kl}i_{lm}i_{mp}\ar[rr] &&  p_{k,lmp}.\\
  }$$

\end{cor}
\dem
Let us consider the following commutative diagram : 
$$
\xymatrix{
\Sigma_{k} \hook[r]&T_{k} \hook[r] & X\\
&T_{k}\cap \Sigma_{l}  \hook[u]^{j_{}} \hook[r]_{i_{l,k}} & \Sigma_{l} \hook[u]_{i_{l}}.
}$$
The composition of the base-change map with the equivalence defined in the proposition \ref{2equivalence} is an equivalence from $i_{kl}$ to $p_{k,l}$. \\
To define the other equivalences we proceed in the same way. \\
The existence of the isomorphisms of functors and their commutations is assured by the fact that the functors $\eta_{l} : i_{km}\rightarrow i_{kl}i_{lm}$ and $\eta_{l} : p_{k,m}\rightarrow p_{k,lm}$ are defined using the $2$-adjunction and by the fact that the base-change map is $2$-functorial.  
\cqfd
\begin{Def}
Let us denote by $\SSS^c_{\Sigma}$ the $2$-category defined as follows.
\begin{itemize}
\item[$\bullet$] The objects are the data  :
\begin{itemize}
\item for every $\Sigma_{k}$, a representation, $\alpha_{k}$, of $\Pi_{2}(\Sigma_{k})$,
\item for every pair $(\Sigma_{k},\Sigma_{l})$ such that $\Sigma_{k}\subset\overline{\Sigma}_{l}$ a functor, $F_{kl}$, :
 $$F_{kl} : \alpha_{k} \longrightarrow p_{k,l}\alpha_{l},$$
\item for every triple $(\Sigma_{k}, \Sigma_{l}, \Sigma_{m})$ such that $\Sigma_{k}\subset \overline{\Sigma}_{l}\subset \overline{\Sigma}_{m}$, an isomorphism of functor $f_{klm}$ visualized by:
$$
\xymatrix @!0 @C=1.3cm @R=0.4cm {\alpha_{k} \ar[rrr]^{F_{lk}}  \ar[ddddd]_{F_{mk}} &&& p_{k,l}(\alpha_{l}) \ar[ddddd]^{p_{k,l}(F_{ml})} \\
~\\
& & \ar@{=>}[ld]_\sim^{f_{klm}}\\
&~\\
\\
p_{k,m}(\alpha_{m}) \ar[rrr]_{i^{-1}_k\eta_{lm}} &&& p_{k,lm}(\alpha_{m})
}$$ such that the following diagram commutes :
$$\xymatrix@!0 @C=3cm @R=0.9cm {  
                       p_{k,lm} (F_{pm}) \circ p_{k,l} (F_{ml} )\circ F_{lk} \ar[rr]^{Id \bullet f_{klm}} \eq[d]^{}&  &  p_{k,lm}(F_{pm}) \circ \eta_l \circ F_{mk} \ar[ddd] \\
                       p_{k,l}\big(p_{l,m}(F_{pm}\big) \circ F_{ml}) \circ F_{lk} \ar[dd]_{(p_{k,l})f_{lmp}\bullet Id_{F_{lk}}}&  & \\
                       &&\\
                       p_{k,l}(\eta_{m} \circ F_{pk}) \circ F_{lk} \eq[d]&&\eta_{l}\circ p_{k,l}(F_{pm})\circ F_{mk}\ar[ddd]^{Id \bullet f_{kmp}} \\
                       p_{k,l}(\eta_{m}) \circ  p_{k,l}(F_{pl}) \circ F_{lk}\ar[dd]_{Id \bullet f_{klp}}  \\
                       \\
                        p_{k,l}(\eta_{m}) \circ \eta_{l} \circ F_{pk} \ar[rr]&  & \eta_l \circ \eta_{m} \circ F_{pk}.
                       }$$   
\end{itemize}
\item[$\bullet$] The $1$-morphisms from $\big(\{\alpha_{k}\}, \{F_{kl}\}, \{f_{klm}\}\big)$ to $\big(\{\alpha'_{k}\}, \{F'_{kl}\}, \{f'_{klm}\}\big)$ are the data  : 
\begin{itemize}
\item for every stratum $\Sigma_{k}$, a functor : $G_{k} : \alpha_{k}\rightarrow \alpha_{k}'$,
\item for every pair $(k,l)$ such that $0\leq k<l\leq n$, an isomorphism of functors~:
$$g_{kl} : F'_{lk}\circ G_{k} \buildrel\sim\over\rightarrow p_{k,l}G_{l}\circ F_{lk}$$
such that the following diagram commutes :
$$\xymatrix@!0 @C=6cm @R=2cm {  
p_{k,l}F'_{ml}\circ F'_{lk}\circ G_{k} \ar[r]^{Id \bullet g_{lk}} \ar[d]_{f_{mlk}\bullet Id} & p_{k,l}F_{ml}'\circ p_{k,l}G_{l}\circ F_{lk}\ar[d]^{p_{k,l}g_{ml} \bullet Id}\\
\eta_{l}\circ F'_{mk}\circ G_{k} \ar[d]_{Id \bullet g_{mk}} & i_{kl}i_{lm}G_{m} \circ p_{k,l}F_{ml}\circ F_{lk} \ar[d]^{Id \bullet f_{mlk}}\\
\eta_{l}\circ p_{k,m}G_{m}\circ F_{mk}\ar[r]&  p_{k,lm}G_{m} \circ \eta_{l}\circ F_{mk} 
}$$
\end{itemize}
\item[$\bullet$] the $2$-morphisms from the $1$-morphism $\big(\{G_{k}\}, \{g_{kl}\}\big)$ to the $1$-morphism $\big(\{G'_{k}\}, \{g_{kl}\}\big)$ are the data for every $0\leq k\leq n$ of a morphism of functors of stacks $\phi_{k} : G_{k}\rightarrow G'_{k}$, such that the following diagram commutes :
$$\xymatrix@!0 @C=6cm @R=2cm { 
F'_{kl}\circ G_{k} \ar[r]^{g_{kl}} \ar[d]_{Id \bullet \phi_{k}} & p_{k,l}G_{l}\circ F_{lk}\ar[d]^{p_{k,l}\phi_{l}\bullet Id}\\
F'_{kl}\circ G'_{k} \ar[r]_{g'_{kl}} & p_{k,l}G_{l}' \circ F_{lk}
}$$
\end{itemize}

\end{Def}
\begin{thm}
The $2$-category, $\SSSt_{\Sigma}^c$ of constructible stacks relatively to $\Sigma$ is equivalent to the $2$-category $\SSS_{\Sigma}^c$.
\end{thm}
\dem
Let us denote by $\SSS$ the image of the $2$-category of constructible stacks relatively to the stratification $\Sigma$ trough $R_{\Sigma}$. 
Using the corollary \ref{cor}, and because of the definition of a constructible stacks relatively to $\Sigma$, it is easy to show that these $2$-category, $\SSS$, is equivalent to the $2$-category whose objects are the data of a family $\big(\{\Lc_{j}\}, \{L_{jk}\}\{l_{jkm}\}\big)$, where $\Lc_{j}$ is a locally constant stack on $\Sigma_{j}$, $L_{jk}$ is a morphism of locally constant stacks $L_{jk}  : \Lc_{j}\rightarrow p_{j,k}\Lc_{k}$ and $l_{jkm}$ is an isomorphism of functors : 
$$
\xymatrix @!0 @C=1.3cm @R=0.4cm {\Lc_{j} \ar[rrr]^{L_{jk}}  \ar[ddddd]_{L_{jm}} &&& p_{j,k}(\alpha_{l}) \ar[ddddd]^{p_{j,k}(F_{mk})} \\
~\\
& & \ar@{=>}[ld]_\sim^{l_{jkm}}\\
&~\\
\\
p_{j,m}(\Lc_{m}) \ar[rrr]_{\eta_{k}} &&& p_{j,km}(\Lc_{m})
}$$ 
satisfying the same commutation  conditions as in the definition of $\SSS_{\Sigma}$
and where the $1$-morphisms the $2$-morphisms are defined in the same way as in the definition of $\SSS_{\Sigma}$.

Now, applying the $2$-monodromy to each locally constant stack, and thank to the theorem \ref{projection}, we show that the $2$-category $\SSS$ is equivalent to $\SSS t^c_{\Sigma}$. 
\cqfd

\nocite{Toen}
\nocite{Macl}
\nocite{Bo}
\nocite{Breen}
\nocite{Gro1}
\bibliographystyle{plain}
\bibliography{biblio}

\end{document}